\documentstyle[12pt]{book}

\makeatletter

\let\ASpecial\special
\newif\if@postscript \@postscripttrue
\def\nopostscript{\let\ASpecial\@ignore\global\@postscriptfalse}

\def\donthyphenate{\pretolerance=10000\hyphenpenalty=10000\exhyphenpenalty=10000}

\let\dropem\relax

\def\aspell#1{\wlog{ASPELL-#1}}

\def\newarray#1{%
\@verbose{NEW COUNTER #1}%
\newcounter{#1@ct}%
\@namedef{#1}##1{\@nameuse{#1@##1}}}

\def\nextelt#1#2{\global\stepcounter{#1@ct}%
\edef\@aktemp{#2}
\global\expandafter\let\csname#1@\arabic{#1@ct}\endcsname\@aktemp}

\def\breakarrayloop{\global\@doarrayloopfalse}

\def\forarray#1{\global\@doarraylooptrue
\bgroup\@verbose{ARRAY FOR LOOP: #1}%
\@namedef{#1}{\@nameuse{#1@\arabic{array}}}
\setcounter{array}{0}%
\@forarray{#1}}

\long\def\@forarray#1#2{\if@doarrayloop\ifnum\c@array<\@nameuse{c@#1@ct}%
\stepcounter{array}%
\@verbose{    Loop step \arabic{array}}%
{#2}
\@forarray{#1}{#2}%
\else
\@verbose{LOOP END}\egroup\fi
\else
\@verbose{Loop Broken}\egroup\fi}

\long\def\fornum#1#2{\setcounter{tempcnta}{#1}%
\setcounter{array}{0}%
\def\aktemp{#2}\@fornum}

\def\@fornum{\ifnum\c@array<\c@tempcnta
\stepcounter{array}%
\@verbose{FORNUM: \arabic{array}}%
\aktemp
\@fornum
\fi}

\def\ps@num@bot{\ps@empty\def\@oddfoot{\hfill\foliofont\thefolio\hfill}\let\@evenfoot\@oddfoot}
\def\footnumbers{\ps@num@bot}

\def\ps@infot{\ps@empty\def\@oddfoot{\hfill\thefolio\hfill\@littleinfobox}\let\@evenfoot\@oddfoot}
\def\ps@infob{\ps@empty\def\@oddhead{\hfill\@littleinfobox}\let\@evenhead\@oddhead}

\def\@testdef #1#2#3{\def\@tempa{#3}\expandafter \ifx \csname #1@#2\endcsname
 \@tempa  \else \@tempswatrue \@verbose{Label ``#2'' has apparently changed}\fi} 

\newdimen\@pimwid
\newdimen\pimgap

\def\putinmargin{\@ifnextchar[{\@putinmargin}{\@putinmargin[]}}
\def\@putinmargin[#1]#2{\ifvmode\leavevmode\fi\vadjust{\noindent\makebox(0,0)[#1r]{\begin{minipage}[b]{\@pimwid}#2\end{minipage}\hskip\pimgap}}}

\newdimen\foliohang 
\newdimen\foliohangcrrct   
\newdimen\fromfolio    

\newdimen\@runheadwid

\def\@evenhead{\hspace*{-\foliohang}\hspace*{-\@sechang}\hbox
to\@runheadwid{{\foliofont
\thefolio}\hskip\fromfolio\outrunheadfill{\markfont\leftmark}\inrunheadfill}}

\def\@oddhead{\hskip-\@sechang\hskip\foliohangcrrct\hbox
to\@runheadwid{\inrunheadfill{\markfont\rightmark}\outrunheadfill\hskip\fromfolio{\foliofont\thefolio}}}

\let\mark@uppercase\relax
\def\uppercasemarks{\let\mark@uppercase\uppercase}

\def\runheadsin{\let\inrunheadfill\hfill
\let\outrunheadfill\relax}

\def\runheadsout{\let\outrunheadfill\hfill
\let\inrunheadfill\relax}

\def\runheadscenter{\let\inrunheadfill\hfill\let\outrunheadfill\hfill}

\runheadsin

\def\thefolio{\colorize{\thepage}}

\def\nopagenumbers{\let\thefolio\relax}

\def\@afterheading{\global\@nobreaktrue
      \everypar{\if@nobreak
                   \global\@nobreakfalse      
                   \clubpenalty \@M
                   \if@afterindent \else {\setbox0=\lastbox}\fi
                 \else \clubpenalty \@clubpenalty
                    \everypar{}\fi
}\@nopgbk[4]\smallmaxskip}

\newdimen\xtrabtwheads
\def\@startsection#1#2#3#4#5#6{\if@noskipsec \leavevmode\fi
   \par \@tempskipa #4\relax
   \@afterindenttrue
   \ifdim \@tempskipa <\z@ \@tempskipa -\@tempskipa \@afterindentfalse\fi
   \if@nobreak \vskip\xtrabtwheads\everypar{}\else
     \addpenalty{\@secpenalty}\addvspace{\@tempskipa}\fi \@ifstar
     {\@ssect{#3}{#4}{#5}{#6}}{\@dblarg{\@sect{#1}{#2}{#3}{#4}{#5}{#6}}}}

\newdimen\secnumdim
\def\@sect#1#2#3#4#5#6[#7]#8{\ifnum #2>\c@secnumdepth
     \let\@svsec\@empty\else
     \refstepcounter{#1}\edef\@svsec{\csname the#1\endcsname\hskip\secnumdim}\fi
     \@tempskipa #5\relax
      \ifdim \@tempskipa>\z@
        \begingroup #6\relax
          \@hangfrom{\hskip #3\relax\@svsec}{\interlinepenalty \@M
\csname#1@uc\endcsname{#8}%
\par}%
        \endgroup
       \csname #1mark\endcsname{#7}\addcontentsline
         {toc}{#1}{\ifnum #2>\c@secnumdepth \else
                      \protect\numberline{\csname the#1\endcsname}\fi
                    #7}\else
        \def\@svsechd{#6\hskip #3\relax  
                   \@svsec \csname#1@uc\endcsname{#8}\csname #1mark\endcsname
                      {#7}\addcontentsline
                           {toc}{#1}{\ifnum #2>\c@secnumdepth \else
                             \protect\numberline{\csname the#1\endcsname}\fi
                       #7}}\fi
     \@xsect{#5}}

\let\section@uc\relax
\let\subsection@uc\relax
\let\subsubsection@uc\relax

\newskip\pre@secskip\newskip\post@secskip

\def\section{\@verbose{SECTION}\ignoreprevskip\dropem
\@startsection{section}{1}{\z@}{\pre@secskip}{\post@secskip}{\sectype\secfont}}

\def\sectype{\noindent\hskip-\@sechang}

\def\sectypea{\noindent$\!$\hskip-\@sechang\setunits{1}%
\real@picture(0,0)
\put(0,-.16){\@loadps{asbar.ps}[5.5,.08]}
\real@endpicture}

\def\sectypec{\noindent}
\def\sectyped{\noindent\hskip-\@sechang\hbox to0pt{\underline{\hskip\fullpgwid}}}

\newdimen\@sechang

\newskip\pre@subsec \newskip\post@subsec

\def\subsection{\@verbose{SUBSECTION}\ignoreprevskip\dropem\@startsection{subsection}{2}{\z@}{\pre@subsec}{\post@subsec}{\subsectype\subsecfont}}

\let\subsectype\sectype

\newskip\pre@subsub \newskip\post@subsub

\def\subsubsection{\@verbose{SUBSUBSECTION}\@startsection{subsubsection}{4}{\z@}{\pre@subsub}{\post@subsub}{\subsubsectype\subsubsecfont}}

\let\subsubsectype\sectypec

\let\stayfigonly\relax
\let\real@minipage\minipage
\let\real@endminipage\endminipage
\def\minipage{\bgroup\stayfigonly\@ifnextchar[{\real@minipage}{\real@minipage[t]}}
\def\endminipage{\real@endminipage\egroup}

\def\caption{\@ifundefined{@captype}{\@caperr}{}\refstepcounter\@captype \@dblarg{\@caption\@captype}}
\def\@caperr{\typeout{ASTeX WARNING: caption USED OUTSIDE OF FLOAT
ENVIRONMENT...USE TCAPTION OR FCAPTION}}

\def\hline{\noalign{\ifnum0=`}\fi\tablecolor \hrule \@height
.5\arrayrulewidth \@black\futurelet\@tempa\@xhline}
\def\vline{\tablecolor{\vrule \@width \arrayrulewidth}\@black}
\def\@arrayrule{\@addtopreamble{\hskip -.5\arrayrulewidth 
   \tablecolor{\vrule \@width .5\arrayrulewidth}\@black\hskip -.5\arrayrulewidth}}
\let\tablecolor\relax

\newdimen\tabletopht \newdimen\tablebotht

\def\@array[#1]#2{\@tempdima\arraystretch\ht\strutbox
\@tempdimb\arraystretch\dp\strutbox
\advance\@tempdima\tabletopht
\advance\@tempdimb\tablebotht
\setbox\@arstrutbox=\hbox{\vrule 
     height\@tempdima
     depth\@tempdimb
     width\z@}\@mkpream{#2}\edef\@preamble{\halign \noexpand\@halignto
\bgroup \tabskip\z@ \@arstrut \@preamble \tabskip\z@ \cr}%
\let\@startpbox\@@startpbox \let\@endpbox\@@endpbox
  \if #1t\vtop \else \if#1b\vbox \else \vcenter \fi\fi
  \bgroup \let\par\relax
  \let\@sharp##\let\protect\relax \lineskip\z@\baselineskip\z@\@preamble}

\def\boxeqnarray{\setbox\@tempboxa=\vtop\bgroup
\stepcounter{equation}\let\@currentlabel=\theequation
\global\@eqnswtrue
\global\@eqcnt\z@\tabskip\@centering\let\\=\@eqncr
\halign to \textwidth\bgroup\@eqnsel\hskip\@centering
  $\displaystyle\tabskip\z@{##}$&\global\@eqcnt\@ne 
  \hskip 2\arraycolsep \hfil${##}$\hfil
  &\global\@eqcnt\tw@ \hskip 2\arraycolsep $\displaystyle\tabskip\z@{##}$\hfil 
   \tabskip\@centering&\llap{##}\tabskip\z@\cr}

\def\endboxeqnarray{\@@eqncr\egroup
\global\advance\c@equation\m@ne\global\@ignoretrue\egroup\box\@tempboxa\par}

\@namedef{boxeqnarray*}{\def\@eqncr{\nonumber\@seqncr}\boxeqnarray}
\@namedef{endboxeqnarray*}{\nonumber\endboxeqnarray}

\def\verbose{\hbadness=10000
\vbadness=10000
\hfuzz36pt
\vfuzz46pt
\def\@verbose##1{\typeout{ASTeX Comment: ##1}}%
\@verbose{VERBOSE MODE SET}}
\def\@verbose#1{}

\def\showfilename{\global\@filenametrue}

\def\put@filename{\@ifundefined{curr@filename}{}{\makebox(0,0)[br]{\fboxsep2pt\fbox{\labfont\curr@filename}\ \ }}}

\newif\if@filename \global\@filenamefalse

\@ifundefined{suffixflag}
	{\def\suffixflag{\relax}}
	{\@verbose{Using suffixflag!}}

\@ifundefined{suffix}
	{\def\suffix{}}
	{\@verbose{Using the suffix \suffix!}}

\def\ignoreprevskip{\ifvmode\@tempdima\lastskip\ifdim\@tempdima>\z@
\@verbose{ignoreprevskip: skip ignored}\vskip-\@tempdima\fi\fi}

\def\maxskip#1{\ifdim\lastskip=\m@gicjuju
\@verbose{MAXSKIP IGNORED}%
\else
\ifdim\lastskip=2\m@gicjuju
\@verbose{Small Maxskip}%
\vskip2pt
\else
\@tempskipa#1\par\advance\@tempskipa-\lastskip
\ifdim\@tempskipa>\z@
\vskip\@tempskipa\fi\fi\fi}
\def\maxspace#1{\@tempskipa#1
\advance\@tempskipa-\lastskip
\ifdim\@tempskipa>\z@\vspace*{\@tempskipa}\fi}

\newdimen\m@gicjuju \m@gicjuju .0834798pt
\def\nomaxskip{\vskip-\m@gicjuju\vskip\m@gicjuju}
\def\smallmaxskip{\vskip-2\m@gicjuju\vskip2\m@gicjuju}

\def\oddnewpage{\@verbose{ODDNEWPAGE Called Here}%
\newpage
\ifodd\c@page
\@verbose{the new page is odd, so I'm leaving it}%
\else
\@verbose{the new page is even...doing another newpage!}%
\thispagestyle{empty}
\phantom{a}\newpage\fi}

\newif\if@savepagenum \@savepagenumfalse
\def\savepagenum{\global\@savepagenumtrue}

\newwrite\@pagefile

\def\autopagenum#1{\def\@lstpg{#1}}		

\newif\if@oddfp \@oddfptrue

\def\startoddpage{\global\@oddfptrue}
\def\startanypage{\global\@oddfpfalse}

\def\okbreak{\penalty -200}
\def\condbreak{\if@nobreak\@verbose{condbreak superceded by NOBREAK}\else
\@verbose{CONDITIONAL BREAK}\penalty-10005\fi}
\def\badbreak{\penalty 90000}

\newdimen\@whitedim
\newbox\@whitebox

\def\whiteout{\@ifstar{\@whiteout}{\@@whiteout}}

\def\@@whiteout#1{{\setbox\@whitebox=\hbox{#1}%
\theprintingcolor
\gray[1]%
\@tempdima\wd\@whitebox \@tempdimb\ht\@whitebox \@tempdimc\dp\@whitebox
\advance\@tempdima2\@whitedim \advance\@tempdimb\@whitedim
\advance\@tempdimc\@whitedim
\hbox{\vrule height \@tempdimb depth\@tempdimc
width\@tempdima
\hskip-\wd\@whitebox \hskip-\@whitedim
\notgray
\@black
\usebox\@whitebox}}}

\def\@whiteout#1{{\@whitedim0pt\@@whiteout{#1}}}

\def\gray{\@ifnextchar[{\@gray}{\@gray[\gray@default]}}
\def\@gray[#1]{\aftergroup\notgray\ASpecial{ps: #1 setgray}}

\def\notgray{\ASpecial{ps: 0 setgray}}

\def\across{\@ifnextchar<{\col@across}{\@across}}

\def\col@across<#1>{\setcounter{across}{0}\@ifnextchar[{\col@ownacross<#1>}{\col@ownacross<#1>[a]}}

\def\col@ownacross<#1>[#2]{\smallskip
\def\@acclab{\csname al@#2\endcsname}%
\@tempdima\linewidth
\divide\@tempdima by #1\relax
\setcounter{qcol}{0}%
\@notfirstfalse
\def\item{\ifnum\c@qcol=#1\hfil\egroup\hfil\egroup
\par\noindent\setcounter{qcol}{0}\fi
\stepcounter{qcol}\relax
\ifnum\c@qcol=1\relax
\noindent
\hbox to\linewidth\bgroup
\@notfirsttrue
\else\hfil\egroup
\fi
\hbox to\@tempdima\bgroup\refstepcounter{across}{\@acclab}}

\def\endacross{\if@notfirst\hfil\egroup\fi\hfil\egroup\par\smallskip}}

\def\@across{\@notfirstfalse\setcounter{across}{0}\@ifnextchar[{\@ownacross}{\@ownacross[a]}}

\def\@ownacross[#1]{%
\def\@acclab{\csname al@#1\endcsname}%
\def\item{\if@notfirst\hskip2pc\fi\@notfirsttrue\refstepcounter{across}\@acclab\nobreak}
\noindent}

\newcounter{across}\def\theacross{\@acclab}
\newif\if@notfirst

\def\al@a{(\alph{across})\penalty\@M\hskip1em}
\def\al@i{(\roman{across})\penalty\@M\hskip1em}
\def\al@I{\Roman{across}.\penalty\@M\hskip1em}

\def\endacross{\par}

\long\def\cfbox#1{\@color\fbox{\@black#1\@color}\@black}

\def\defbox{\maxskip{\intextsep}\@ifnextchar[{\@@defbox}{\@@defbox[\textwidth]}}

\long\def\@@defbox[#1]#2{{\advance\fboxsep-\fboxrule
\newbox\@defbox
\global\@tempdima#1\global\advance\@tempdima-2\fboxsep%
\setbox\@defbox=\hbox{\begin{minipage}[b]{\@tempdima}
#2
\end{minipage}}
\global\@tempdimb\ht\@defbox \global\advance\@tempdimb\fboxsep
\par
\noindent
\centerline{\hbox{\@defboxshade\cfbox{\hfill\box\@defbox\hfill}}}
\par
\vskip\intextsep}}

\def\defboxshade{.8}

\def\@defboxshade{\@color\makebox(0,0)[bl]{\ASpecial{ps: \defboxshade\space setgray}%
\@tempdima\wd\@defbox\advance\@tempdima2\fboxsep\advance\@tempdima2\fboxrule
\vrule height\@tempdimb width\@tempdima depth\fboxsep
\ASpecial{ps: 0 setgray}}\@black}

\newdimen\ruleboxside 
\def\rulebox{\maxskip{\intextsep}\@ifstar{\@findrulebox}
{\@ifnextchar [ {\@rulebox}{\@rulebox[\cutpgwid]}}}

\long\def\@rulebox[#1]#2{\okbreak
\@tempdima#1 \advance\@tempdima-2\fboxsep \advance\@tempdima-2\fboxrule
\badbreak\par
{\centering
\par
\cfbox{\begin{minipage}{\@tempdima}
\hskip\ruleboxside
\@tempdima\textwidth \advance\@tempdima-2\ruleboxside
\begin{minipage}{\@tempdima} 
#2
\end{minipage}
\end{minipage}}
\par}
\vskip\intextsep}

\long\def\@findrulebox#1{\@tempdima\cutpgwid \advance\@tempdima-2\fboxsep
\setbox\@tempboxa=\vbox{\begin{minipage}[b]{\@tempdima}#1\end{minipage}}%
\@tempdimb\ht\@tempboxa
\ifdim\@tempdimb>\baselineskip
\@tryless{#1}\relax
\else
\setbox\@tempboxa=\hbox{#1}\relax
\@tempdima\wd\@tempboxa \advance\@tempdima 2\fboxsep
\typeout{THIS IS NOT AN ERROR....THIS IS NOT AN ERROR.
This is the output of rulebox* which attempts to find
the best width of a piece of text.  This one turns out
to fit on one line...its width is }%
\showthe\@tempdima
\rulebox[\@tempdima]{#1}%
\fi}

\long\def\@tryless#1{\advance\@tempdima-\fboxsep
\setbox\@tempboxa=\vbox{\begin{minipage}[b]{\@tempdima}#1\end{minipage}}%
\ifdim\ht\@tempboxa>\@tempdimb
\advance\@tempdima 4\fboxsep
\typeout{THIS IS NOT AN ERROR....THIS IS NOT AN ERROR.
This is the output of rulebox* which seeks to find the best width for
a given set of text.  Please hard put this value next time.
I have calculated the width of this text as }%
\showthe\@tempdima
\rulebox[\@tempdima]{#1}%
\else
\@tryless{#1}\fi}

\def\fullpagesize{\par
\noindent\hskip-\@sechang\hskip-\@totalleftmargin
\minipage[b]{\fullpgwid}\unskip\unskip
\fps@centering\par}

\let\fps@centering\relax

\def\pullback{\noindent\hskip-\@sechang\hfil}

\let\footnotecolor\relax
\def\footnoterule{\footnotecolor{\kern-3\p@
  \hrule width 8pc \kern 2.6\p@}\@black} 

\def\breakhere{\@verbose{BREAKHERE}%
\par\ignoreprevskip\break}

\def\pareject{{\parfillskip0pt\eject}\noindent}

\def\runover{\typeout{ASTeX Warning: RUNOVER CALLED}%
\@runover}
\def\@runover{\global\pagetotal-\pagegoal}

\long\def\ignore#1{\typeout{ASTeX WARNING: IGNORE IN PROGRESS}}
\long\def\@ignore#1{}

\def\newpagination#1{\@verbose{NEW PAGINATION: #1}%
\@namedef{during#1pag}{\@ignore}%
\@namedef{#1breakhere}{\relax}%
\@namedef{#1pareject}{\relax}%
\@namedef{#1runover}{\relax}%
\@namedef{#1ignore}{\relax}%
\@namedef{#1vspace}{\@ifstar{\@ignore}{\@ignore}}%
\@namedef{#1pagination}{\@namedef{#1breakhere}{\breakhere}%
\@namedef{#1pareject}{\pareject}%
\@namedef{during#1pag}{\relax}%
\@namedef{#1runover}{\runover}%
\@namedef{#1vspace}{\vspace}%
\@namedef{#1ignore}{\ignore}}}

\def\docalculations{\@verbose{Calculating Page Format}%
\@maxdepth\maxdepth
\@@parindent\parindent
  \@sechang\fullpgwid \advance\@sechang-\cutpgwid
\ifdim\@sechang<\z@
\typeout{ASTeX Format Error: fullpgwid is less than cutpgwid;}%
\typeout{ASTeX Format Error: this may produce poor results.}%
\fi
\@pimwid\@sechang \advance\@pimwid-\pimgap
 \oddsidemargin\paperwd \advance\oddsidemargin-\fullpgwid
 \divide\oddsidemargin by2 \@othermargin\oddsidemargin
 \advance\oddsidemargin\@sechang
 \evensidemargin\oddsidemargin
 \textwidth\cutpgwid
\@runheadwid\fullpgwid \advance\@runheadwid\foliohang
\@totaltop\topmargin\advance\@totaltop\headheight
\advance\@totaltop\headsep}

\newif\if@asformat \@asformatfalse

\newdimen\@@parindent 

\def\asfonts#1{\def\@asfonts{#1}}

\def\fmtfile{as-format}
\def\fntfile{as-fonts}
\def\defsfile{as-defs}

\def\asformat#1{\edef\holdcode{\the\catcode`@}\catcode`@=11
\typeout{}
\typeout{ASTeX Format #1 Selected}
\def\@aktemp{#1}\def\@@aktemp{DEFAULT}%
\ifx\@aktemp\@@aktemp
\fmt@def
\else
\input{\fmtfile#1}%
\fi
\docalculations
\@ifundefined{@asfonts}{\def\@asfonts{DEFAULT}\typeout{ASTeX WARNING: No
font file specified...using DEFAULT}}{}
\def\@@aktemp{DEFAULT}
\ifx\@@aktemp\@asfonts
\typeout{Using default fonts!}
\fnt@def
\else
\input{\fntfile\@asfonts}
\fi
\if@postscript
\ASpecial{! userdict begin
/ASscale { UnitScale div scaletopta mul } def
/curveweight { \ps@specdef{curveweight}{.75} ASscale } def
/axisweight  { \ps@specdef{axisweight}{.25} ASscale } def
/gridweight  { \ps@specdef{gridweight}{.15} ASscale } def
/contourweight { \ps@specdef{contourweight}{1} ASscale} def
/curveweight3d { \ps@specdef{curveweightTd}{1} ASscale } def
/axisweight3d  { \ps@specdef{axisweightTd}{.5} ASscale } def
/gridweight3d  { \ps@specdef{gridweightTd}{.3} ASscale } def
/lightshade  { \ps@specdef{lightshade}{.9} } def
/medshade    { \ps@specdef{medshade}{.75} } def
/darkshade   { \ps@specdef{darkshade}{.6}  } def
/lightline { \ps@specdef{lightline}{.6} } def
/medline { \ps@specdef{medline}{.4} } def
/darkline { \ps@specdef{darkline}{0} } def
/trap { \ps@specdef{trap}{.35} } def 
/calcdashlen { \ps@specdef{calcdashlen}{4} ASscale } def
/calcgaplen  { \ps@specdef{calcgaplen}{2} ASscale } def
/ticklength  { \ps@specdef{ticklength}{3} ASscale } def
/calcpointsize   { \ps@specdef{calcpointsize}{4} ASscale } def
/regsize { \ps@specdef{regmarksize}{.5} } def
/showgridbig { \ps@specdef{showgridbig}{.5} } def 
/showgridsmall { \ps@specdef{showgridsmall}{.75} } def
/ticker { ticklength Mathabs div add } def
/tickerer { /AKtemp exch def 
            /AAKtemp { AKtemp ticklength Mathabs div sub } def 
            0 AAKtemp ge { AKtemp }{ AAKtemp } ifelse } def
/ASarrowlen { \ps@specdef{ASarrowlen}{2} } def
/ASarrowwid { \ps@specdef{ASarrowwid}{1} } def
/ASarrowdict 22 dict def
ASarrowdict begin
/ASmtrx matrix def
end
/ASarrow
{ASarrowdict begin
   /ASarrowscale exch def
   /tipy exch def /tipx exch def
   /taily exch def /tailx exch def
   /dx tipx tailx sub def
   /dy tipy taily sub def
   /arrowlength dx dx mul dy dy mul add
     sqrt def
   /headlength 0.04 ASarrowlen mul ASarrowscale mul UnitScale div def 
   /base arrowlength headlength sub def
   /halfthickness currentlinewidth .5 mul def
   /halfheadthickness 0.02 ASarrowwid mul ASarrowscale mul UnitScale div def
   arrowlength 0.0 eq {  
        tailx taily halfthickness 0 360 arc
        }
   {/savematrix ASmtrx currentmatrix def
   /angle dy dx atan def
   tailx taily translate
   angle rotate
   base 0 ge { 0 halfthickness neg moveto
   base halfthickness neg lineto
   base halfthickness lineto
   0 halfthickness lineto
   closepath fill } if
   base halfheadthickness neg moveto
   arrowlength 0 lineto
   base halfheadthickness lineto
   closepath fill
   savematrix setmatrix} ifelse
   end
   } def
/ASarrowcent
{ASarrowdict begin
   /ASarrowscale exch def
   /diry exch def /dirx exch def
   /centy exch def /centx exch def
   /dx dirx centx sub def
   /dy diry centy sub def
   /dirlength dx dx mul dy dy mul add
     sqrt def
   /headlength 0.04 ASarrowlen mul ASarrowscale mul UnitScale div def 
   /halfheadthickness 0.02 ASarrowwid mul ASarrowscale mul UnitScale div def
   dirlength 0.0 eq {  
        centx centy halfthickness 0 360 arc
        }
   {/savematrix ASmtrx currentmatrix def
   /angle dy dx atan def
   centx centy translate
   angle rotate
   headlength .5 mul neg halfheadthickness neg moveto
   headlength .5 mul 0 lineto
   headlength .5 mul neg halfheadthickness lineto
   closepath fill
   savematrix setmatrix} ifelse
   end
   } def
/UnitScale { 1 } def
/lineshade { setgray } bind def
/calcrotate { /AKangle exch def 
currentpoint /AKtempY exch def
/AKtempX exch def
gsave
AKtempX AKtempY translate
AKangle -1 mul rotate 
AKtempX -1 mul AKtempY -1 mul translate 
} def
end}

\@input{psfig.sty}\@ifundefined{psfig}{\@input{psfig}}{}
\@ifundefined{psfig}{\typeout{ASTeX Warning: COULD NOT FIND
PSFIG...the command loadps expects to use this package.}}{%
\@ifundefined{@noisyfalse}{}{\@noisyfalse}}
\catcode`@=11

\def\@color{\ASpecial{ps: ascolor}}
\def\@black{\ASpecial{ps: asblack}}
\def\theprintingcolor{\ASpecial{ps: printingcolor}}

\def\setcolorrestore{\ASpecial{ps: setcolorrestore}}
\def\colorrestore{\ASpecial{ps: colorrestore}}
\def\@setfigcolor{\ASpecial{ps: picinit}}
\def\@resetcolor{\ASpecial{ps: resetcolor}}
\def\blackonly{\ASpecial{header=\psdir/ASTeX/blackset.ps}}
\def\coloronly{\ASpecial{header=\psdir/ASTeX/colorset.ps}}
\long\def\colorize##1{\@color##1\@black}

\def\color{\@color\aftergroup\@black}

\ASpecial{! userdict begin
/ascolor{
   doingsep {
	colormemory {
		coloralready not {
			switchtocolor
		} if
		/coloralready { true } def
		/blackalready { false } def
	}{
		switchtocolor
	} ifelse
   } if
   } def
/switchtocolor{
   kotrap{
	blackonly{
		knockout
	} if
	coloronly{
		1 currentgray eq{
			0 setgray
		} if
	} if
    }{
		colorshift
    } ifelse
  } def
/asblack{
   doingsep {
	colormemory {
		blackalready not {
			switchtoblack
		} if
		/coloralready { false } def
		/blackalready { true } def
	}{
		switchtoblack
	} ifelse
   } if
   } def
/switchtoblack{
   kotrap{
	coloronly{
		knockout
	} if
	blackonly{
		1 currentgray eq{
			0 setgray
		} if
	} if
    }{
		blackshift
    } ifelse
  } def
/colorshift{
	cstrans shiftfactor mul 
	cstrans shiftfactor mul 
	translate
 } def
/blackshift{
	-1 cstrans mul shiftfactor mul
	-1 cstrans mul shiftfactor mul
	translate
 } def
	/knockout{
		 1 setgray
		 /AKtemp{ currentlinewidth scaletoptb mul UnitScale mul
		 trap 2 mul 
		 sub } def
		AKtemp 0 ge { AKtemp scaletoptb div UnitScale div }{ 0 } ifelse
		setlinewidth
	 } def
 /printingcolor { doingsep {
			  coloronly { TeXcolor } if
			  blackonly { TeXblack } if
		 } if
                 } def
/setcolorrestore { doingsep {
			coloralready { /colorrestore {TeXcolor} def } if
			blackalready { /colorrestore {TeXblack} def } if 
			} if
		} def
/colorrestore { } def
 /pco	 { 	  /blackalready { false } def
		  /coloralready { true } def
		  /coloronly    { true } def
		  /blackonly    { false } def } def
 /pbo	 { /blackalready { true } def
		  /coloralready { false } def
		  /blackonly    { true } def
		  /coloronly    { false } def } def
 /picinit { printingcolor } def
 /resetcolor {  coloronly {
				/coloralready { true } def
                                /blackalready { false } def } if 
                blackonly { 
				/coloralready { false } def
                                /blackalready { true } def } if } def
 /printcoloronly { /bop-hook { pco } def 
			/doingsep{ true } def
			/cstrans { translateby } bind def } def
 /printblackonly { /bop-hook { pbo } def 
			/doingsep{ true } def
			/cstrans { translateby } bind def } def
 /blackalready { true } def
 /coloralready { true } def
 /blackonly { false } def
 /coloronly { false } def
 /cstrans { 0 } def
 /translateby { 10000 } def
 /doingsep { false } def
 /scaletopta { 1 } def
 /scaletoptb { 1 } def
 /swlcmd{setlinewidth}def
 /shiftfactor { 1 }def
 /kotrap { false } def
 /putcommandshere {} def
 putcommandshere
end}

\def\PSswitches##1##2##3##4##5##6{\@verbose{NEW COLOR COMMAND: ##1}%
\ASpecial{! userdict begin
/##1switches { /kotrap {##2} def
	/colormemory {##3} def
	/scaletopta{##4} def
	/scaletoptb{##5} def
	/shiftfactor{##6} def } def
end}}

\PSswitches{TeX}{false}{true}{1}{1}{1}
\PSswitches{M}{true}{false}{1}{Mathabs}{.05}
\PSswitches{F}{true}{true}{1.1}{1.1}{1}

\ASpecial{! userdict begin TeXswitches end}

\else
\def\color{}
\def\@color{}
\def\@black{}
\def\theprintingcolor{}
\def\setcolorrestore{}
\def\colorrestore{}
\def\@setfigcolor{}
\def\@resetcolor{}
\def\blackonly{}
\def\coloronly{}
\def\colorize{}
\fi
\@input{\defsfile#1}%
\suffixflag
\@asformattrue\catcode`@=\holdcode}

\def\nocolor{\def\@nocolor{TRUE}}
\def\nographics{\def\@nographics{TRUE}}

\let\real@document\document
\def\document{\endgroup\@verbose{BEGIN DOCUMENT}%
\make@pmboxes
\if@asformat\@verbose{ASFormat has been specified}%
\else\typeout{ASTeX: No Format has been
specified...assuming DEFAULTS}\asformat{DEFAULT}\fi
\if@apbi\@ifundefined{issuenumdir}{\typeout{ASTeX Error:  You have
specified ``auto page by issue'' (so that each issue begins on the
page after the previous issue.)  However, to use this you must specify
a directory \string\issuenumdir.  This macro has not been defined.}}{}%
\fi
\if@apbc\@ifundefined{homedir}{\typeout{ASTeX Error: You have spcified
``auto page by chapter'' (so that each chapter begins on the page
after the previous one).  However, to use this, you must specify a
home directory \string\homedir.  This macro has not been defined.}}{}\fi
\begingroup\real@document}

\long\def\Ncolumns#1#2#3{%
\@verbose{Ncolumns #1}%
\ifvmode\ifdim\pagegoal>\paperht\@verbose{Ncolumns: PAGE
EMPTY...adding hardtop}\hardtop\par\vskip-2\baselineskip\else
\vskip0pt plus 12pt\fi
\else
\ifdim\pagegoal>\z@
\vskip0pt plus 12pt\fi\fi
\par\@tempdima\pagegoal \advance\@tempdima-\pagetotal
\setbox\@tempboxa=\vbox{\vspace*{0pt}\par#3\par\vspace*{0pt}}%
\advance\@tempdima-\ht\@tempboxa
\ifdim\@tempdima<.5in\@verbose{NCOLUMN forcing page
break}\pagebreak[4]\par\hardtop\par\vskip-2\baselineskip\fi
\@Ncolht\textheight
\@verbose{Ncolumns: text height is \the\@Ncolht}%
\advance\@Ncolht-\ht\@tempboxa
\@verbose{Ncolumns: Leaving room for heading gives \the\@Ncolht}%
\advance\@Ncolht-\pagetotal
\@verbose{Ncolumns: and leaving room for what is already on the page: \the\@Ncolht}%
\@tempdima\pagetotal
\advance\@tempdima-2\pagegoal
\global\pagetotal\@tempdima
{
\box\@tempboxa\nopagebreak[4]\par\nopagebreak[4]%
\@Ncolwid\fullpgwid\advance\@Ncolwid\@Ncolsep
\advance\@Ncolwid-#1\@Ncolsep
\divide\@Ncolwid by #1
\hsize\@Ncolwid \linewidth\hsize
\global\setbox\@Ncolbox=\vbox{\pretolerance=5000\vspace*{0pt}\par#2}}%
\@Ncolsplit\@Ncolbox{#1}%
\ifdim\pagetotal<\z@
\@tempdima\pagetotal
\advance\@tempdima2\pagegoal
\pagetotal\@tempdima\fi}

\def\@Ncolsplit#1#2{\par
\setcounter{tempcnta}{0}%
\noindent\hskip-\@sechang\vbox{\hbox to\fullpgwid{\splitmaxdepth0pt\boxmaxdepth0pt\@doNsplits{#1}{#2}}\hardtop\par}\par
\@Ncolht\textheight
\ifdim\ht#1>\z@
\pagebreak[4]%
\@Ncolsplit#1{#2}
\fi
\par\okbreak}
\newif\if@shrt \@shrtfalse
\newif\if@allnow \@allnowfalse

\def\@doNsplits#1#2{\ifnum\c@tempcnta<#2\stepcounter{tempcnta}%
\ifnum\c@tempcnta=1%
\@tempdima\ht#1\divide\@tempdima by #2\advance\@tempdima10pt%
\ifdim\@tempdima<\@Ncolht\@Ncolht\@tempdima
\global\@shrttrue
\@verbose{Ncolumn: short columns}\fi
\fi
\if@shrt\ifnum\c@tempcnta=#2\global\@allnowtrue\fi\fi
\@verbose{Ncolumns -- column \arabic{tempcnta}}%
\if@allnow
\global\@allnowfalse\global\@shrtfalse
\setbox\@tempboxa=\vbox to \@Ncolht{\unvbox#1\vfill}\@emptyit#1%
\else
\setbox\@tempboxa=\vsplit#1 to \@Ncolht%
\setbox\@tempboxa=\vbox{\unvbox\@tempboxa}%
\fi
\@verbose{column height \the\ht\@tempboxa\ and goal \the\@Ncolht}%
\ifdim\ht\@tempboxa>\z@\ifdim\wd\@tempboxa=\z@
\addtocounter{tempcnta}{-1}\else
\setbox\@tempboxa=\vbox to \@Ncolht{\unvbox\@tempboxa}%
\ifnum\badness>\@M
\@verbose{TOO BIG TO FIT!}%
\setbox\@tempboxa=\vbox{\unvbox\@tempboxa}%
\fi
\hbox to \@Ncolwid{\vtop{\hardtop\par\@labbox\@tempboxa}}
\ifnum\c@tempcnta<#2\hskip\@Ncolsep\fi
\fi\fi
\@doNsplits{#1}{#2}\fi}

\def\@addsomemore#1{\advance\@tempdima1pc%
\setbox\@tempboxb=\vsplit#1 to \@tempdima
\@verbose{Addsomemore: this column is too short!}%
\ifdim\ht\@tempboxb>\z@
\setbox\@tempboxa=\vbox{\unvbox\@tempboxa\par
\noindent
\unvbox\@tempboxb}%
\else
\@verbose{Addsomemore: need more}%
\fi}

\def\indexhere{\@ifnextchar*{\star@ih}{\nostar@ih}}
\def\star@ih*{\let\indexhead\relax\def\@indname{book}\as@indexhere}
\def\nostar@ih{\def\@indname{\jobname}\as@indexhere}

\def\as@indexhere{\def\theindex{\par\parskip\z@\relax\let\item\@idxitem}
\def\endtheindex{\par}
\def\indexspace{\par \vskip 10pt plus 1pt\relax}
\@ifnextchar[{\@indexhere}{\@indexhere[3]}}

\def\@indexhere[#1]{\Ncolumns{#1}{\@input{\@indname.ind}}{\indexhead}}

\def\indexhead{\section{Index}\vspace*{\post@secskip}}

\long\def\extratoc#1{
\if@toc
\addcontentsline{toc}{extra}{#1}
\else
\typeout{ASTeX Warning: The
command `\string\extratoc' was called, but no toc file is open}
\fi}

\def\numberline#1{#1\ \ }

\newif\if@tcs \global\@tcstrue

\def\contentsline#1#2#3{\@verbose{CONTENTS: #1 #2}%
\par{\@ifundefined{tcf@#1}{\secfont}{\csname
tcf@#1\endcsname}\csname
tcl@#1\endcsname{#2}{#3}}\par
\if@tcs
\@ifundefined{tcs@#1}{\smallskip}{\csname tcs@#1\endcsname}
\else
\global\@tcstrue
\fi}

\def\tcl@a#1#2{\sectypea #1 \hfill#2}
\def\tcl@b#1#2{\noindent\hskip-\@sechang #1 \hfill#2}
\def\tcl@c#1#2{\sectypec #1 \hfill#2}
\def\tcl@d#1#2{\sectyped #1 \hfill#2}
\def\tcl@e#1#2{\noindent\hskip-\@sechang #1 \dotfill#2}
\def\tcl@f#1#2{\sectypec #1 \dotfill#2}
\def\tcl@g#1#2{\noindent\llap{#2\ \ }#1}
\def\tcl@h#1#2{\noindent#1\hfill}
\def\tcl@yst#1#2{\sectype #1 \hfill#2}
\def\tcl@ysst#1#2{\subsectype #1 \hfill#2}
\def\tcl@ignore#1#2{\global\@tcsfalse}
\def\tcl@center#1#2{\begin{center}#1\end{center}}

\def\tcf@article{\normalsize}

\def\authtitla#1#2{{\tcf@auth #1}\\{\tcf@titl#2}}
\def\authtitlb#1#2{{\tcf@titl #2}\\{\tcf@auth #1}}
\def\authtitlc#1#2{{\tcf@auth #1 -- }{\tcf@titl#2}}
\def\authtitld#1#2{{\tcf@titl #2 -- }{\tcf@auth #1}}

\let\authtitl\authtitlb

\let\tcf@auth\it \let\tcf@titl\bf
\let\tcl@article\tcl@f
\let\tcl@chapter\tcl@yst
\let\tcl@section\tcl@ysst
\let\tcl@subsection\tcl@ignore
\let\tcl@subsubsection\tcl@ignore
\let\tcl@extra\tcl@center

\def\tochead{\title{Table of Contents}}

\def\tableofcontents{\@ifnextchar*{\@alltoc}{{\newpage\let\thepage\thetocpage\tochead\@filetoc{\jobname}\maketoc\nostretch\newpage}%
\if@apbi\setpagenumber\else\if@apbc\else
\oddnewpage
\setcounter{page}{1}\fi\fi}}

\def\@alltoc*{{\let\thepage\thetocpage\tochead\forarray{chap}{\@filetoc{\homedir/\chap\suffix}}}\nostretch\newpage}

\def\thetocpage{\roman{page}}

\def\@filetoc#1{\@input{#1.toc}}

\def\maketoc{\@starttoctoc\global\@toctrue\@verbose{MAKE TABLEOFCONTENTS}
\toc@doc\global\let\maketoc\relax}

\def\toc@doc{\global\let\toc@enddocument\enddocument\global\def\enddocument{\@opentoc\toc@enddocument}}

\def\@starttoctoc{\newwrite\tf@toc}

\def\@opentoc{\@verbose{Creating \jobname.toc}\immediate\openout\tf@toc=\jobname.toc\relax}

\def\nameasfont{\@ifnextchar*{\@@nameasfont}{\@nameasfont}}

\def\@@nameasfont*#1#2{\@namedef{#1}{\@fontnote{#1}\asfont@init#2}}

\def\@fontnote#1{}
\def\asfontnotify{\def\@fontnote##1{\typeout{ASTeX FONT: ##1}}}
\newskip\@tempbls

\@ifundefined{reset@font}{\let\asfont@init\relax}{\let\asfont@init\reset@font}

\def\@nameasfont#1#2#3{\@verbose{naming font: #1}%
\expandafter\font\csname @#1\endcsname=#2 at #3pt%
\@namedef{#1}{\@fontnote{#1}\asfont@init\protect\@nameuse{mathsize@#1}\protect\@nameuse{fontextra@#1}%
\baselineskip\@nameuse{bls@#1}\protect\@nameuse{@#1}}%
\@ifnextchar({\@verbose{   Baselineskip Specified}\def@bls{#1}{#3}}
{\@tempdima#3pt\divide\@tempdima by 5%
\@tempdimb#3pt
\advance\@tempdimb\@tempdima\def@bls{#1}{#3}(\the\@tempdimb)}}

\def\def@bls#1#2(#3){\expandafter\newdimen\csname bls@#1\endcsname
\csname bls@#1\endcsname#3
\@ifnextchar[{\@verbose{   Math Size Specified}\@mathsize{#1}}%
{\@verbose{   Determining Best Math Size}%
\@tempdimb100pt%
\setbox\holdpicbox=\hbox{\csname @#1\endcsname x}%
\forarray{fontsizes}{%
\setbox\@tempboxa=\hbox{\expandafter\csname\fontsizes\endcsname x}
\@tempdima\ht\holdpicbox \advance\@tempdima-\ht\@tempboxa
\ifdim\@tempdima<\z@\breakarrayloop
\else
\ifdim\@tempdima<\@tempdimb\setcounter{tempcnta}{\c@array}
\@tempdimb\@tempdima
\fi
\fi
}
\found@mathsize{#1}}}

\def\found@mathsize#1{
\ifcase \c@tempcnta
\typeout{Please edit definition of found@mathsize}
\or
\@namedef{mathsize@#1}{\tiny}
\or
\@namedef{mathsize@#1}{\scriptsize}
\or
\@namedef{mathsize@#1}{\footnotesize}
\or
\@namedef{mathsize@#1}{\small}
\or
\@namedef{mathsize@#1}{\normalsize}
\or
\@namedef{mathsize@#1}{\large}
\or
\@namedef{mathsize@#1}{\Large}
\or
\@namedef{mathsize@#1}{\LARGE}
\or
\@namedef{mathsize@#1}{\huge}
\or
\@namedef{mathsize@#1}{\Huge}
\or \typeout{Please edit definition of found@mathsize}\fi
\@ifnextchar<{\@fontextra{#1}}{\@fontextra{#1}<\relax>}}

\def\@mathsize#1[#2]{%
\@verbose{Setting #1 mathsize to #2}%
\@tempswatrue
\edef\@@aktemp{#2}%
\forarray{fontsizes}{\edef\@aktemp{\fontsizes}%

\ifx\@aktemp\@@aktemp\global\@tempswafalse\breakarrayloop\fi}
\if@tempswa
\typeout{ASTeX Format Warning: The font declaration of #1 specifies a
math size which is not a standard LaTeX size}
\fi
\@namedef{mathsize@#1}{\csname #2\endcsname}
\@ifnextchar<{\@fontextra{#1}}{\@fontextra{#1}<\relax>}}

\def\@fontextra#1<#2>{\@namedef{fontextra@#1}{#2}}

\newdimen\paperwd \newdimen\@othermargin
\newdimen\@totaltop
\newdimen\@@parindent
\newdimen\fullpgwid \newdimen\cutpgwid
\newdimen\paperht
\newarray{fontsizes}
	\nextelt{fontsizes}{tiny}
	\nextelt{fontsizes}{scriptsize}
	\nextelt{fontsizes}{footnotesize}
	\nextelt{fontsizes}{small}
	\nextelt{fontsizes}{normalsize}
	\nextelt{fontsizes}{large}
	\nextelt{fontsizes}{Large}
	\nextelt{fontsizes}{LARGE}
	\nextelt{fontsizes}{huge}
	\nextelt{fontsizes}{Huge}
\newarray{so}
\newbox\@Ncolbox
\newcounter{array}
\newcounter{tempcnta}
\newcounter{tempcntb}
\newdimen\@Ncolht
\newdimen\@Ncolsep \@Ncolsep 1pc
\newdimen\@Ncolwid
\newdimen\@tempdima
\newdimen\@tempdimd
\newdimen\mlinebot
\newdimen\mlinetop
\newif\if@doarrayloop 
\newif\if@ins \global\@insfalse
\newif\if@toc \@tocfalse
\newif\if@oddfp \@oddfptrue

\newif\if@qinput \@qinputtrue
\newarray{exlist}
\newwrite\@exfile
\newif\if@makesolpic \@makesolpictrue

\newdimen\qcolwid
\newcounter{qcol}

\newarray{chap}
\newarray{chapfile}
\newarray{chapdir}

\def\newexerhead#1#2{\@namedef{exerhead@#1}{#2}}
\newskip\@solskip \newskip\@exampskip 
\def\chapheada{\thispagestyle{empty}\noindent%
\hskip-\@sechang\begin{fullpagesize}
\centering
\if@chapapp
{\@color\chapfont {CHAPTER~\numberword{chapter}}\@black}
\else
{\@color\chapfont {APPENDIX~\thechapter}\@black}
\fi
\vskip 10pt

\unitlength\textwidth
\@loadps{asbar.ps}[1,.04]
\vskip 31pt

\colorize{
\begin{minipage}{25pc}
\centering
\Huge\chaptitlfont \thetitle
\end{minipage}
}
\vskip4.5pc

\centerline{\@loadps{asbar.ps}[.9,.02]}
\vskip 8pt
\begin{minipage}{.8\textwidth}
\parindent\@@parindent
\noindent
\thesummary
\end{minipage}
\vskip 1pc
\centerline{\@loadps{asbar.ps}[.9,.02]}
\end{fullpagesize}
\vfill
\begin{fullpagesize}
\centering
\unitlength 1pt
\real@picture(0,0)
\put(0,-18){\makebox(0,0)[t]{\colorize{\foliofont \thepage}}}
\real@endpicture
\end{fullpagesize}
\newpage}

\def\chapheadc{\rulebox{%
\centering

\chapfont \if@chapapp Chapter\else Appendix\fi\ \thechapter}

\rulebox{\centering

\chaptitlfont
\thetitle}

\rulebox{
\thesummary}

\thispagestyle{num@bot}\medskip\medskip\medskip\vspace*{0pt}\newpage}

\def\chapheadd{\thispagestyle{num@bot}{\flushright\noindent\underline{\hbox
to\textwidth{\hfill\chapfont \if@chapapp Chapter \else Appendix \fi \thechapter}}\par\smallskip
\hfill\begin{minipage}{.85\textwidth}\donthyphenate
\flushright
\chaptitlfont \thetitle
\end{minipage}}\par\bigskip}

\def\chapheade{\vspace*{-2\baselineskip}\par\noindent\unitlength\oddsidemargin\real@picture(0,0)%
\put(-1,0){\raise\@totaltop\hbox{\makebox(0,0)[tl]{\vrule
height \paperht width \paperwd}}}
\@tempdima\paperht\advance\@tempdima-\@totaltop
\put(0,0){\gray[1]\vrule height \@totaltop depth \@tempdima width 2pt}
\unitlength\textheight
\@tempdima\cutpgwid\advance\@tempdima\@othermargin
\put(0,-.3){\gray[1]\vrule height 2pt width \@tempdima}
\put(0,-.2){\makebox(0,0)[bl]{\hskip1pc\begin{minipage}{.7\textwidth}
\donthyphenate
\flushleft
\par\noindent\whiteout{\chapfont \if@chapapp CHAPTER 
\numberword{chapter}
\else
APPENDIX \thechapter \fi}\end{minipage}}}
\put(0,-.4){\makebox(0,0)[tl]{\hskip1pc\begin{minipage}{.7\textwidth}%
\donthyphenate
\flushleft
\gray[1]%
\noindent{\chaptitlfont \thetitle\par}\notgray
\end{minipage}}}\real@endpicture\vspace*{2pt}\newpage}

\newif\if@chapapp \@chapapptrue

\def\appendices{\@verbose{APPENDICES instead of CHAPTERS}%
\global\@chapappfalse\global\def\thechapter{\Alph{chapter}}}

\def\issuenum#1{\def\theissue{#1}\@ifundefined{lastpagedir}{}{\savepagenum
\@verbose{Issue Number Set -- Will seek last page of previous issue.}}}

\def\thejournal{THE NAME OF YOUR JOURNAL}

\def\infobox{\small\bf \thejournal \\
\small Vol. \thevolume\ No. \theissue\ \thepubdate\\
\small Pages \thepage -- \thelastpage}

\def\putinfobox[#1#2]{\useinfobox
\def\ps@firstpage{\csname
ps@info#1\endcsname}\def\@infotb{#1}%
\@makeinfolr{#2}}

\def\@infobox{\hbox to \textwidth{\@leftinfofill 
\begin{tabular}{c}\infobox\end{tabular}\@rightinfofill}\global\@infoboxherefalse}

\def\useinfobox{\global\def\article@infobox{
\vspace*{0pt}\nointerlineskip\par
\@verbose{Making Infobox}%
\@tempdima\pagegoal
\setbox\@tempboxa=\vbox{\let\pageref\relax\let\ref\relax\@infobox}
\advance\@tempdima-\ht\@tempboxa
\advance\@tempdima-\infospace
\global\pagegoal\@tempdima\par
\thispagestyle{firstpage}
\global\@infoboxheretrue}}

\def\article@infobox{\thispagestyle{firstpage}}

\def\@makeinfolr#1{\def\@@aktemp{#1}%
\let\@leftinfofill\hfill
\let\@rightinfofill\hfill
\def\@aktemp{l}
\ifx\@aktemp\@@aktemp
\let\@leftinfofill\relax
\fi
\def\@aktemp{r}
\ifx\@aktemp\@@aktemp
\let\@rightinfofill\relax
\fi}

\def\autopage#1{\@input{#1\suffix.lstpg}%
\@ifundefined{@lstpg}%
	{\typeout{ASTeX Warning: Last Page Reference Not Found!!!
		 Check #1\suffix.}\typeout{ASTeX Warning: Starting at page 1.}%
\setcounter{page}{1}}
	{\setcounter{page}{\@lstpg}%
	\if@oddfp
	\ifodd\c@page
		\@verbose{the first page is odd}%
	\else
		\addtocounter{page}{1}%
		\@verbose{skipping an even page so that first page
will be odd}%
	\fi\fi}%
\def\autopage##1{\typeout{ASTeX Warning: AUTOPAGE called more than once}}}

\def\as@label[#1]#2{\def\@aktemp{-}\def\@@aktemp{#2}%
\ifx\@aktemp\@@aktemp
\@verbose{AS@LABEL request...no label here}%
\else
\@bsphack\if@filesw {\let\thepage\relax
   \def\protect{\noexpand\noexpand\noexpand}%
\xdef\@gtempa{\write\@auxout{\string
      \asnewlabel{#1\as@nlextra}{#2}{{\@currentlabel}{\thepage}}}}}\@gtempa
   \if@nobreak \ifvmode\nobreak\fi\fi\fi\@esphack\fi}

\def\as@nlextra{}

\def\asnewlabel#1{\newlabel}

\newif\if@cap
\newif\if@fcap
\newbox\tfinputbox
\newbox\tfcapbox
\newdimen\tcapskip \newdimen\fcapskip 
\def\finput{\global\@fcaptrue\global\def\as@labtag{finput}\tfinput}
\def\tinput{\global\@fcapfalse\global\def\as@labtag{tinput}\tfinput}

\def\tfinput{\@ifstar{\nocap@tfinput}{\cap@tfinput}}

\def\cap@tfinput{\@ifnextchar[{\global\def\as@labtag{no-file}\let\curr@filename\relax
\lab@tfinput}{\nolab@tfinput}}

\def\nolab@tfinput#1{%
\def\curr@filename{#1}%
\lab@tfinput[#1]{\input{#1}}}

\long\def\lab@tfinput[#1]#2{%
\def\@aktemp{-}\def\@@aktemp{#1}%
\ifx\@aktemp\@@aktemp\else
\def\curr@filename{#1}\fi
\setbox\tfinputbox=\hbox{\if@filename\put@filename\fi#2}%
\capwd@tfinput[#1]}

\def\capwd@tfinput[#1]{%
\@ifnextchar<{\cw@tfinput[#1]}{\cw@tfinput[#1]<\wd\tfinputbox>}}

\long\def\cw@tfinput[#1]<#2>#3{%
\setbox\tfcapbox=%
\hbox{\begin{minipage}{#2}
\if@fcap\fcaption{#3}\else\tcaption{#3}\fi\as@label[\as@labtag]{#1}\end{minipage}}\end@tfinput}

\def\nocap@tfinput{\@ifnextchar[{\lnc@tfinputa}{\nlnc@tfinput}}

\def\nlnc@tfinput#1{\lnc@tfinputa[#1]{\input{#1}}}

\long\def\lnc@tfinputa[#1]#2{\def\curr@filename{#1}\setbox\tfinputbox=\hbox{\if@filename\put@filename\fi#2}\lnc@tfinput[#1]}

\def\lnc@tfinput[#1]{\setbox\tfcapbox=\hbox{\begin{minipage}{\wd\tfinputbox}%
\if@fcap\fnumber\else\tnumber\fi\as@label[\as@labtag]{#1}\end{minipage}}\end@tfinput}

\def\end@tfinput{\if@fcap\makefinput\else\maketinput\fi}

\def\makefinput{\ifdim\wd\tfinputbox>\wd\tfcapbox
\@tempdima\wd\tfinputbox
\else
\@tempdima\wd\tfcapbox
\fi
\begin{minipage}{\@tempdima}\centering
\par
\centerline{\@finputtop}\par\vskip\fcapskip\centerline{\@finputbot}\par
\vspace*{0pt}\end{minipage}}

\def\figontop{\def\@finputtop{\box\tfinputbox}\def\@finputbot{\box\tfcapbox}}
\def\figonbot{\def\@finputbot{\box\tfinputbox}\def\@finputtop{\box\tfcapbox}}

\def\maketinput{\ifdim\wd\tfinputbox>\wd\tfcapbox
\@tempdima\wd\tfinputbox
\else
\@tempdima\wd\tfcapbox
\fi
\begin{minipage}{\@tempdima}\centering
\par
\centerline{\@tinputtop}\par\vskip\tcapskip\centerline{\@tinputbot}\par
\vspace*{0pt}\end{minipage}}

\def\tabontop{\def\@tinputtop{\box\tfinputbox}\def\@tinputbot{\box\tfcapbox}}
\def\tabonbot{\def\@tinputbot{\box\tfinputbox}\def\@tinputtop{\box\tfcapbox}}

\figontop
\tabonbot
\let\fcaptype\centering
\let\fcapfont\small

\def\post@fcaptype{:~}

\long\def\fcaption#1{\refstepcounter{figure}\vbox{\fcaptype\@fnumber\post@fcaptype\fcapfont#1\par}}
\def\fnumber{\refstepcounter{figure}\vbox{\fcaptype\par\@fnumber\par}}
\def\@fnumber{{\figlegfont Figure~\thefigure}}

\let\tcaptype\centering
\let\tcapfont\small

\def\post@tcaptype{:~}
\long\def\tcaption#1{\refstepcounter{table}\vbox{\tcaptype\@tnumber\post@tcaptype\tcapfont#1\par}}
\def\tnumber{\refstepcounter{table}\vbox{\tcaptype\par\@tnumber\par}}
\def\@tnumber{{\tablegfont Table~\thetable}}

\def\bibheada{\noindent\hskip-\@sechang{\subsecfont
References}\small}

\let\bibhead\bibheada

\let\bibtype\bibtypea

\def\@frpb#1{{\parfillskip\z@\hskip0pt plus 1fil\nolinebreak\hfill\hspace*{.5in}\llap{#1}\par}}

\def\proofboxa{\@frpb{\rule{7pt}{7pt}\hskip-8pt\raisebox{1pt}{\hbox{\gray\rule{7pt}{7pt}\notgray}}}}
\def\proofboxb{{\@frpb{\fboxsep\z@\fbox{\gray\rule{9pt}{9pt}\notgray}}}}
\def\proofboxc{\@frpb{\fboxrule.5pt\fboxsep\z@\fbox{\phantom{Ii}}}}
\def\proofboxd{{\@frpb{\gray\rule{9pt}{9pt}\notgray}}}
\def\proofboxe{\nopagebreak[4]\leavevmode
\nopagebreak[4]\vadjust{\hfill\real@picture(0,0)
\setlength{\unitlength}{1pc}
\gray
\put(0,-.4){\thicklines\colorize{\line(-1,0){14.9}}}
\notgray
\real@endpicture}}
\let\proofbox\proofboxc
\def\claimtypea#1{\noindent\hskip-\@sechang{\hbox{\claimheadfont#1}}\claimbodyfont}

\def\claimtyped#1{\noindent\hskip-\@sechang{\colorize{\hbox{\claimheadfont#1}}\claimbodyfont}}
\def\claimtypee#1{\noindent
\hskip-\@sechang
\if@needproof
\unitlength 1pc%
\real@picture(0,0)
\gray
\put(0,1){\thicklines\colorize{\line(1,0){14.9}}}
\notgray
\real@endpicture\fi
{\colorize{\hbox{\claimheadfont#1}}\claimbodyfont}}
\let\claimtype\claimtypea

\def\abstracttypeb{\small\hskip\@@parindent\bf}

\newskip\preabstractskip
\newskip\postabstractskip
\let\abstracttype\abstracttypa
\newdimen\infospace \infospace.5pc

\def\setpagenumber{\if@apbi
	\newpage
	\setcounter{tempcnta}{\theissue}%
	\setcounter{tempcntb}{\c@page}%
	\addtocounter{tempcnta}{-1}%
	\ifnum\c@tempcnta>0%
		\@verbose{Getting Last Page Number from Issue \arabic{tempcnta}}%
		\autopage{\issuenumdir/\ispref\arabic{tempcnta}}
		\ifnum\c@tempcntb>1%
			\addtocounter{tempcntb}{-\c@page}%
			\ifodd\c@tempcntb
				\@verbose{Adding a Blank Page to Correct Parity}%
				\addtocounter{page}{-1}%
				\phantom{a}
				\thispagestyle{empty}%
				\par\break
			\else
				\@verbose{Parity is Correct}%
			\fi
	\fi
	\else
		\@verbose{First Issue...so starting at page one.}
	\fi
\else
	\typeout{ASTeX Warning: Auto Page by ISSUE not Specified}%
\fi}

\def\sectionmark#1{}

\newif\if@apbi \@apbifalse 

\def\articledefaults{\apbi\maketoc
\def\thefigure{\arabic{section}.\arabic{figure}}
\def\thetable{\arabic{table}}
\def\thesection{\arabic{section}}
\def\thesubsection{\arabic{section}.\arabic{subsection}}
\def\theequation{\arabic{equation}}
\let\tcl@chapter\tcl@ignore
\let\tcl@section\tcl@ignore
\let\tcl@subsection\tcl@ignore}

\def\apbi{\global\@apbitrue\global\@apbcfalse\savepagenum}

\def\noautopage{\global\@apbifalse\global\@apbcfalse\global\@savepagenumfalse}

\def\@littleinfobox{\makebox(0,0)[r]{\begin{tabular}{c}\littleinfobox\end{tabular}}}

\let\ps@firstpage\ps@num@bot

\def\littleinfobox{}

\def\@artname{\jobname}

\newif\if@oddartfp  \@oddartfpfalse

\def\startartodd{\global\@oddartfptrue}
\def\startartany{\global\@oddartfpfalse}

\long\def\article#1#2{\@ifnextchar[
{\@article{#1}{#2}}{\@@article{#1}{#2}}}

\long\def\@@article#1#2#3{\@article{#1}{#2}[#3]{#3}} 

\long\def\@article#1#2[#3]#4{\@secinput[#1]
{\do@article{#1}{#2}{#4}{#3}}
{\@verbose{SKIPPING ARTICLE #1}}
}

\long\def\do@article#1#2#3#4{{\let\label\@artlabel
\let\ref\@artref
\let\@bibref\art@bibref
\let\pageref\@artpageref
\if@oddartfp
\oddnewpage
\else
\newpage
\fi
\article@infobox
\if@filename
\def\curr@filename{#1}\put@filename\fi
\markboth{\mark@uppercase{#2}}{\mark@uppercase{#4}}%
\addcontentsline{toc}{article}{\protect\authtitl{#2}{#3}}
\def\@artname{#1}%
\typeout{}%
\typeout{ARTICLE ---> #1}%
\setcounter{section}{0}%
\setcounter{equation}{0}%
\setcounter{figure}{0}%
\setcounter{table}{0}%
\setcounter{footnote}{0}%
\aspell{#1}\input{#1}\@svpg\nostretch\newpage}}

\def\@artlabel#1{\@bsphack\if@filesw {\let\thepage\relax
   \def\protect{\noexpand\noexpand\noexpand}%
   \edef\@tempa{\write\@auxout{\string
      \newlabel{\@artname-#1}{{\@currentlabel}{\thepage}}}}%
   \expandafter}\@tempa
   \if@nobreak \ifvmode\nobreak\fi\fi\fi\@esphack}

\def\@artref#1{\@ifundefined{r@\@artname-#1}{
\@ifundefined{r@#1}{{\bf ??}}{\real@ref{#1}}}{\edef\@tempa{\@nameuse{r@\@artname-#1}}\expandafter
\@car\@tempa \@nil\null}}

\def\@artpageref#1{\@ifundefined{r@\@artname-#1}{
\@ifundefined{r@#1}{{\bf ??}\typeout{REF WARNING: Reference `#1' on
page \thepage \space in the article ``\@artname''
undefined}}{\real@ref{#1}}}{\edef\@tempa{\@nameuse{r@\@artname-#1}}\expandafter
\@cdr\@tempa \@nil\null}}
 
\let\real@ref\ref

\def\@svpg{\label{last page number}}
\def\thelastpage{\pageref{last page number}}

\def\@jnlinfo{\jnlinfo\if@showpages\ pp. \thepage-\pageref{last page number}\fi}

\newif\if@showpages \@showpagesfalse
\def\showpages{\@showpagestrue}

\def\bibliography{\par\dropem\bigskip\bibhead\nopagebreak[4]\enumerate}

\def\@@makebib{\hfill\bibitemtype}
\def\bibitemtype{\bibtype}
\long\def\bibitem[#1]#2{\refstepcounter{enumi}\item[\@@makebib{\arabic{enumi}}]\label{br@#1}#2}
\def\bibref[#1]{\bibtype{\bibeval[#1--@]}}

\def\bibeval[#1--#2]{\@bibeval[#1,@]\def\@aktempo{#2}\def\@aktempi{@}\ifx\@aktempo\@aktempi\else--\bibeval[#2]\fi}

\def\eatthespace#1{\@eatthespace#1\endeatthespace}
\def\@eatthespace{\@ifnextchar*{\@@eatthespace}{\@@eatthespace}}
\def\@@eatthespace#1\endeatthespace{\def\@aktempo{#1}}

\def\@bibeval[#1,#2]{%
\split@bibref[#1|@]%
\def\@aktempo{#2}\def\@aktempi{@}\ifx\@aktempo\@aktempi\else\bibseps\@bibeval[#2]\fi}

\def\bibseps{,~}

\def\split@bibref[#1|#2]{\eatthespace{#1}\@bibref{\@aktempo}\def\@aktempo{#2}\def\@aktempi{@}%
\ifx\@aktempo\@aktempi
\else
\@subbibref[#2]\fi}

\def\@subbibref[#1|@]{\subbibref{#1}}

\def\subbibref#1{~\hbox{\rm (#1)}}

\def\@bibref#1{\@ifundefined{r@br@#1}{??}{\ref{br@#1}}}

\def\art@bibref#1{\@ifundefined{r@\@artname-br@#1}{??}{\ref{br@#1}}}

\def\exactrefs{\@verbose{EXACT REFERENCES
CALLED}\def\@bibref##1{##1}\def\bibitem[##1]{\item[\bibtype{##1}]}}

\def\end@claim{\ignoreprevskip
\if@proofbox
\proofbox\par\maxskip{\@exampskip}%
\else
\global\@proofboxtrue\fi}

\newif\if@proofbox \@proofboxtrue
\def\noproofbox{\global\@proofboxfalse}

\let\endproof\end@claim

\def\endprooffill{\ifvmode\else\linebreak[2]\fi\hskip0pt plus 10 fil}

\newcounter{claim}[section]
\def\theclaim{\thesection.\arabic{claim}}

\newif\if@needproof \@needprooffalse

\def\newclaim#1{\@namedef{state#1}{\state<#1>}%
\@namedef{#1}{\@claim<#1>}\@namedef{end#1}{\end@claim}\@namedef{#1ct}{\claimtype}}

\def\state<#1>{\@ifnextchar[{\@state<#1>}{\spec@state<#1>}}

\def\spec@state<#1>{\@ifnextchar<{\@spec@state<#1>}{\@state<#1>[-]}}

\long\def\@spec@state<#1><#2>#3{\maxskip{\@exampskip}
\makeclaim{#1}{#2}{#3}%
\nopagebreak[4]
\if@needproof\@proof\fi\maxskip{\@exampskip}\global\@needprooffalse}

\long\def\makeclaim#1#2#3{\par\def\@aktemp{#2}\def\@@aktemp{-}%
{\csname #1ct\endcsname{\firstup[#1]\ifx\@aktemp\@@aktemp\else~#2\fi:~~}%
#3\par}}

\def\firstup[#1#2]{\uppercase{#1}#2}

\def\@proof{\global\@nobreakfalse \global\@noskipsectrue
       \everypar{\if@noskipsec \global\@noskipsecfalse
                   \clubpenalty\@M \hskip -\parindent
                   \proof\let\@proofhead\@@proofhead
                    \else \clubpenalty \@clubpenalty
                    \everypar{}\fi}\ignorespaces}
\def\proof{\noindent{\bf \@proofhead:~~}}
\def\@@proofhead{Proof}
\let\@proofhead\@@proofhead
\def\proofhead#1{\def\@proofhead{#1}}

\def\claimtypea#1{\noindent\hskip-\@sechang{\hbox{\claimheadfont#1}}\claimbodyfont}

\def\claimtyped#1{\noindent\hskip-\@sechang\colorize{\hbox{\claimheadfont#1}}\claimbodyfont}
\def\claimtypee#1{\noindent
\hskip-\@sechang
\if@needproof
\unitlength 1pc%
\real@picture(0,0)
\gray
\put(0,1){\thicklines\colorize{\line(1,0){14.9}}}
\notgray
\real@endpicture\fi
\colorize{\hbox{\claimheadfont#1}}\claimbodyfont}

\let\claimtype\claimtypea

\long\def\@state<#1>[#2]#3{\refstepcounter{claim}%
\def\@aktempo{#2}\def\@aktempi{-}%
\ifx\@aktempo\@aktempi\else
\label{#2}%
\fi
\if@needproof\else\maxskip{\@exampskip}\fi
\makeclaim{#1}{\theclaim}{#3}%
\if@needproof\@proof\fi\maxskip{\@exampskip}\global\@needprooffalse}

\def\@claim{\@needprooftrue\maxskip{\@exampskip}%
\state}

\def\note{\@ifnextchar[{\@note}{\@note[Note]}}
\def\@note[#1]{\par\smallskip\noindent{\bf#1:} }

\newclaim{claim}
\newclaim{theorem}
\newclaim{corollary}
\newclaim{lemma}
\newclaim{proposition}

\def\tageqnarray[#1]{\def\@lbltmp{#1}$$\parbox[c]{.9\textwidth}\bgroup
\begin{eqnarray*}}
\def\endtageqnarray{\end{eqnarray*}\egroup\tag\@lbltmp$$}

\def\cite#1{\bibref[#1]}

\long\def\abstract#1{\vskip\preabstractskip\par
\begin{quote}\abstracttype#1\end{quote}
\vskip\postabstractskip\par}

\def\abstracttypeb{\small\hskip\@@parindent\bf}

\newskip\preabstractskip
\newskip\postabstractskip

\let\abstracttype\abstracttypea

\long\def\title#1{{\begin{center}\donthyphenate\chapfont{#1}
\end{center}\par}}

\def\author#1{\smallskip\par\begin{center}{\authorfont #1}\@author}
\def\@author{\@ifnextchar[{\@@author}{\end@author}}
\def\@@author[#1]{\\{\addressfont#1}\@author}
\def\end@author{\end{center}\par\smallskip}

\let\apn@document\document
\def\document{
\apn@document
\@bsphack
\ifnum\c@chap@ct>0%
\edef\@@aktemp{\jobname}
\@tempswafalse
\forarray{chapfile}{\def\@chaptemp{\chapfile\suffix}{\escapechar=-1\xdef\@aktemp{\expandafter\string\csname\@chaptemp\endcsname}}%
\@verbose{\@aktemp = \@@aktemp}%
\ifx\@aktemp\@@aktemp\global\@tempswatrue
\breakarrayloop
\@verbose{Found this file in chaplist}\fi}
\if@tempswa
\maketoc
\addtocounter{array}{-1}%
\if@apbc 
\@verbose{Auto Page by Chapter}%
\savepagenum
\ifnum\c@array>0
\autopage{\homedir/\expandafter\chap{\arabic{array}}}
\else
\@verbose{It is first in the list...starting on page 1}
\fi\fi
\@verbose{Collecting GLOBAL ref files}%
\forarray{chap}{\@input{\homedir/\chap\suffix.grf}}
\fi\fi
\@esphack}

\def\textbookdefaults{\def\sectionmark##1{\markright{\mark@uppercase{\thesection~##1}}}%
\apbc
\let\tcl@article\tcl@ignore
}

\newif\if@apbc \@apbcfalse

\def\apbc{\global\@apbctrue\global\@apbifalse}

\def\collectindex{\newwrite\idx\immediate\openout\idx=book.idx
\def\indexentry##1##2{\immediate\write\idx{\string\indexentry{##1}{##2}}}%
\forarray{chap}{\@input{\homedir/\chap\suffix.idx}}}

\def\@secinput[#1]#2#3{\def\@aststamp{#1}
\ifnum\c@just@sec@ct>0\def\@aktemp{#1}\@tempswafalse
\forarray{just@sec}{\edef\@@aktemp{\just@sec}\ifx\@aktemp\@@aktemp\global\@tempswatrue\breakarrayloop\fi}%
\if@tempswa
{%
\def\c@just@sec@ct{0}
#2}\else{#3}\fi
\else
#2
\fi}

\newarray{just@sec}

\def\secinput{\@ifnextchar*{\star@secinput}{\real@secinput}}

\def\real@secinput#1{\@ifnextchar[{\@real@secinput{#1}}{\one@real@secinput{#1}}}

\long\def\one@real@secinput#1#2{\@real@secinput{#1}[#2]{#2}}

\long\def\@real@secinput#1[#2]#3{\@secinput[#1]{\section[#2]{#3}\aspell{#1}\as@label[secinput]{#1}\input{#1}}{\stepcounter{section}}}

\def\star@secinput*#1{\aspell{#1}\@secinput[#1]{\input{#1}}{}}

\def\fmt@def{%
\articledefaults
\noautopage
\asfonts{DEFAULT}
\putinfobox[b]
\fromfolio 1.5pc
\pre@secskip -12.9pt plus -2.5pt minus -2pt
\post@secskip 8pt
\pre@subsec -10pt plus -2.5pt minus -2pt
\post@subsec 3pt
\pre@subsub 10pt plus 2.5pt minus 2pt
\post@subsub 3pt
\fussy 
\maxdepth=2.2pt
\@maxdepth\maxdepth
\hbadness10000
\vbadness10000
\vfuzz 150pt
\hfuzz 150pt
\baselineskip 12pt
\fboxrule 1pt
\fboxsep 1pc
\tabcolsep 9pt
\tabletopht 2pt
\tablebotht 1pt
\parsep 4.5pt
\parindent 18pt
\@beginparpenalty=2000
\listparindent\parindent
\headsep 2pc
\intextsep 18pt plus 2.5pt minus 13pt %
\parskip\z@ 
\@solskip 22pt plus 2.5pt \@exampskip 10pt plus 2.5pt minus 4pt
\labelsep 1em 
\leftmargini 2.5em 
\def\@listI{\leftmargin\leftmargini \parsep\z@ \topsep\z@ 
\@beginparpenalty10000
\partopsep\z@ 
\itemsep 3pt}
\let\@listi\@listI
\leftmargin\z@
\def\@listii{\leftmargin\leftmarginii
 \labelwidth\leftmarginii\advance\labelwidth-\labelsep
 \topsep 3pt
 \@beginparpenalty50
 \partopsep\z@
 \itemsep\z@
 \parsep\z@}
\arrayrulewidth 1pt
\fcapskip.5pc
\tcapskip.5pc
\setcounter{secnumdepth}{2}
\fullpgwid 5.5in
\cutpgwid 5.25in
\paperwd 8.5in
\topmargin .75in
\textheight 8.25in
\foliohang 0pt
\foliohangcrrct 0pt 
\voffset -1in
\hoffset -1in
\newdimen\@paperheight
\paperht11in
\@whitedim 3pt
\def\gray@default{.8}
\nopostscript
}

\def\fnt@def{%
\nameasfont*{labfont}{\small}
\nameasfont*{chapfont}{\Huge\bf}
\nameasfont*{secfont}{\bf}
\nameasfont*{subsecfont}{\it}
\nameasfont*{subsubsecfont}{\it}
\nameasfont*{chaptitlfont}{\Huge\bf}
\nameasfont*{examphead}{\it}
\nameasfont*{solhead}{\it}
\nameasfont*{foliofont}{\bf}
\nameasfont*{figlegfont}{\small\bf}
\nameasfont*{authorfont}{\bf}
\nameasfont*{addressfont}{\it}
\nameasfont*{claimheadfont}{\bf}
\nameasfont*{claimbodyfont}{\it}
\nameasfont*{markfont}{\footnotesize}
}

\newbox\@pmbox
\newbox\odd@pmbox
\newbox\even@pmbox
\newdimen\dim@regmark

\def\make@pmboxes{%
\global\setbox\odd@pmbox=\hbox{\@textbox\unitlength\oddsidemargin\begin{picture}(0,0)%
\put(-1,0){\raise\@totaltop\hbox{\@paperbox\makebox(0,0)[br]{\@tlcrop}\makebox(0,0)[l]{\hspace*{.5\paperwd}\raise\dim@regmark\hbox{\makebox(0,0)[c]{\@regmark}}}%
\raise-\paperht\hbox{\makebox(0,0)[tr]{\@blcrop}\makebox(0,0)[l]{\hspace*{.5\paperwd}\raise-\dim@regmark\hbox{\makebox(0,0)[c]{\@regmark}}}}%
\raise-.5\paperht\hbox{%
\makebox(0,0)[l]{\hspace*{\paperwd}\hspace*{\dim@regmark}\makebox(0,0)[c]{\@regmark}}%
}}}%
{\unitlength\cutpgwid\advance\unitlength\@othermargin
\put(1,0){\raise\@totaltop\hbox{\makebox(0,0)[bl]{\@trcrop}%
\raise-\paperht\hbox{\makebox(0,0)[tl]{\@brcrop}}}}}
\end{picture}}
\global\setbox\even@pmbox=\hbox{\@textbox\unitlength\oddsidemargin\begin{picture}(0,0)%
\put(-1,0){\raise\@totaltop\hbox{\@paperbox\makebox(0,0)[br]{\@tlcrop}\makebox(0,0)[l]{\hspace*{.5\paperwd}\raise\dim@regmark\hbox{\makebox(0,0)[c]{\@regmark}}}%
\raise-\paperht\hbox{\makebox(0,0)[tr]{\@blcrop}\makebox(0,0)[l]{\hspace*{.5\paperwd}\raise-\dim@regmark\hbox{\makebox(0,0)[c]{\@regmark}}}}%
\raise-.5\paperht\hbox{%
\hspace*{-\dim@regmark}\makebox(0,0)[c]{\@regmark}}}}
{\unitlength\cutpgwid\advance\unitlength\@othermargin
\put(1,0){\raise\@totaltop\hbox{\makebox(0,0)[bl]{\@trcrop}%
\raise-\paperht\hbox{\makebox(0,0)[tl]{\@brcrop}}}}}
\end{picture}}}

\newcounter{sheet}

\def\@makepmbox{\stepcounter{sheet}%
\setbox\@pmbox=\hbox{\ifodd\c@page\usebox\odd@pmbox
\else\usebox\even@pmbox\fi
\makebox(0,0)[bl]{\raise\@totaltop\hbox{\raise\@cropsize\hbox{\asfont@init\topstamp}}}}}

\def\@aststamp{\jobname}

\def\topstamp{\hbox{\rm\asfont@init{\bf Job: \jobname}\ \ {\it Sheet: \arabic{sheet}\
Page: \thepage\/}\ \ (\today)\ \ \@ifundefined{@aststamp}{}{[\@aststamp]}}}

\def\@@tlcrop{{\unitlength\@cropsize\begin{picture}(1,1)%
\put(1,.1){\line(0,1){.9}}
\put(0,0){\line(1,0){.9}}\end{picture}}}

\def\@@trcrop{{\unitlength\@cropsize\begin{picture}(1,1)%
\put(0,.1){\line(0,1){.9}}
\put(.1,0){\line(1,0){.9}}\end{picture}}}

\def\@@blcrop{{\unitlength\@cropsize\begin{picture}(1,1)%
\put(1,0){\line(0,1){.9}}
\put(0,1){\line(1,0){.9}}\end{picture}}}

\def\@@brcrop{{\unitlength\@cropsize\begin{picture}(1,1)%
\put(0,0){\line(0,1){.9}}
\put(.1,1){\line(1,0){.9}}\end{picture}}}

\def\@@regmark{\ASpecial{ps: gsave
currentpoint translate
.5 setlinewidth
0 -72 regsize mul moveto
0 72 regsize mul lineto
stroke
-72 regsize mul 0 moveto
72 regsize mul 0 lineto
stroke
0 0 .8 72 mul regsize mul 0 360 arc
stroke
grestore}}

\newdimen\@cropsize \@cropsize.5in

\def\real@textbox{\makebox(0,0)[tl]{\gray\hskip-\@sechang
\unitlength\textheight
\begin{picture}(0,0)\put(0,0){\line(0,-1){1}}\end{picture}}%
\makebox(0,0)[tl]{\gray\fboxsep\z@\fboxrule.5pt%
\@tempdima\textwidth\advance\@tempdima-\fboxrule
\hskip-.5\fboxrule\fbox{\phantom{\vrule height \textheight width \@tempdima}}}}

\def\real@paperbox{\makebox(0,0)[tl]{\gray\fboxsep\z@\fboxrule.5pt%
\@tempdima\paperwd\advance\@tempdima-\fboxrule
\@tempdimb\paperht\advance\@tempdimb-\fboxrule
\fbox{\phantom{\vrule height \@tempdimb width \@tempdima}}}}

\def\showtextbox{\let\@textbox\real@textbox}
\def\notextbox{\let\@textbox\relax}
\notextbox

\def\showpaperbox{\let\@paperbox\real@paperbox}
\def\nopaperbox{\let\@paperbox\relax}
\nopaperbox

\def\showregmarks{\let\@regmark\@@regmark}
\def\noregmarks{\let\@regmark\relax}
\noregmarks

\def\nostamp{\let\topstamp\relax}

\def\nocropmarks{\let\@tlcrop\relax
\let\@trcrop\relax
\let\@blcrop\relax
\let\@brcrop\relax}

\def\showcropmarks{\let\@tlcrop\@@tlcrop
\let\@trcrop\@@trcrop
\let\@blcrop\@@blcrop
\let\@brcrop\@@brcrop}
\nocropmarks
\nostamp

\def\deg{\hbox{deg }}

\makeatother

\newclaim{proposition}
\def\eqref#1{$(\ref{#1})$}

\def\d{\partial}
\def\C{{\blackboard C}}

\def\four#1{{#1}^{\flat}}
\def\BD{{\blackboard D}}
\def\BDbeta{\BD_{0}}

\makeatletter
\let\c@subsection\c@claim

\let\twelveBbb\bf

\def\blackboard#1{\hbox{\twelveBbb #1}}
\makeatother
\begin{document}

\thispagestyle{num@bot}
\hbox to \textwidth{\hfill\begin{minipage}{2.5in}
\small University of Georgia\\
Mathematics Preprint Series\\ No. 2, Volume 4 (1996)\\ To appear in
{\it Physica D}
\end{minipage}}

\bigskip

\title{Bispectral Darboux Transformations: The
Generalized Airy Case}

\author{Alex Kasman and
Mitchell Rothstein}[Department of Mathematics] [University of
Georgia][Athens, GA 30602]

\abstract{This paper considers Darboux transformations of
a bispectral operator which preserve its bispectrality.  A sufficient
condition for this to occur is given, and applied to the case of
generalized Airy operators of arbitrary order $r>1$.  As a result, the
bispectrality of a large family of algebras of rank $r$ is
demonstrated.  An involution on these algebras is exhibited which
exchanges the role of spatial and spectral parameters, generalizing
Wilson's rank one bispectral involution.  Spectral geometry and the
relationship to the Sato grassmannian are discussed.}

\section{Introduction}

An ordinary differential operator, $L(x,\d_x)$, is said to be
bispectral \cite{firstbisp} if it has an eigenfunction $\psi(x,z)$
satisfying a pair of eigenvalue equations
\begin{equation} L(x,\d) \psi(x,z)=f(z) \psi(x,z)\qquad Q(z,\d_z)
\psi(x,z)=g(x)
\psi(x,z)\label{bispeqn}
\end{equation} for non-constant functions $f$ and $g$. This property
was investigated by Duistermaat and Gr\"unbaum \cite{DG}, who
identified all bispectral Schr\"odinger operators: $L=\d^2+V(x)$. In
\cite{W}, Wilson considered the more general question of classifying
bispectral commutative rings of ordinary differential operators.  That
is, he considered rings ${\cal L}$ and ${\cal Q}$ of operators in the
variable $x$ and $z$ respectively such that there exists a function
$\psi(x,z)$ for which Equation~\eqref{bispeqn} holds for all
$L\in{\cal L}$ and $Q\in{\cal Q}$.

Following \cite{BC}, one defines the {\it rank\/} of a commutative
ring of ordinary differential operators to be the greatest common
divisor of the orders of its elements. Wilson classified all maximal
bispectral rings of rank one:

\statetheorem{{\bf[Wilson,\cite{W}]} Let ${\cal R}$ be a maximal rank one
commutative algebra of ordinary differential operators. Then ${\cal
R}$ is bispectral if and only if its spectral curve is a rational
curve with no singularities other than cusps.}

This result was proved, in part, through the introduction of an
involution on a subset $Gr^{ad}$ of the grassmannian $Gr^1$ generally
used in the investigation of the KP hierarchy \cite{SW}.  By the work
of Sato, Segal, Wilson and others \cite{SW} one knows how to associate
to a point $W\in Gr^1$ a commutative ring of operators and the Baker
function $\psi_W(x,z)$ which is a common eigenfunction of the
operators in the ring. Then $Gr^{ad}$ is the subspace of $Gr^1$
consisting of suitably normalized points whose corresponding ring is
the affine coordinate ring of a rational curve with cusps, minus a
point at $\infty$. Wilson showed that $Gr^{ad}$ is equipped with an
involution, $\beta:Gr^{ad}\to Gr^{ad}$, with the property that
\begin{equation}
\psi_{\beta(W)}(x,z)=\psi_W(z,x)
\end{equation} and consequently that the corresponding rings are
bispectral. (Essentially, this is the involution which exchanges the
rings ${\cal L}$ and ${\cal Q}$.)

The classification of bispectral algebras in higher rank appears to be
more subtle \bibref[Bakalov,Gr2].  There is, however, a well-known
method for constructing new commutative algebras of differential
operators from a given one: {\it Darboux transformation}.  (See, for
example, \cite{Wright,Z1,ZM} for previous results relating these
transformations to the bispectral problem.)  Under suitable
hypotheses, described below, Darboux transformations preserve
bispectrality.  In this paper we will explore this method in the case
of generalized Airy operators.

The general approach is the following.  Let ${\cal D}$ denote the ring
of ordinary differential operators with coefficients in $\C(x)$.
Given a function $\psi(x,z)$, consider the equation
\begin{equation}\label{pair}
T(x,\d)\psi=\four{T}(z,\d_z)\psi\ ,
\end{equation}
for $T$ and $\four{T}$ in ${\cal D}$.  Define
\begin{equation}
{\cal A}_{\psi}\subset{\cal D}
\end{equation} to be the set of operators $T$ for which there exists
$\four{T}$ such that Equation~\eqref{pair} holds.  Note that ${\cal
A}_{\psi}$ is a subalgebra of ${\cal D}$.  If $\psi$ is sufficiently
general, in particular if $\psi$ is a bispectral eigenfunction
satisfying equation
\eqref{bispeqn}, then the map
\begin{equation} T\mapsto\four{T}
\end{equation} is an anti-isomorphism from ${\cal A}_{\psi}$ to ${\cal
A}_{\tilde\psi}$, where
\begin{equation}
\tilde\psi(x,z)=\psi(z,x)\ .
\end{equation}
Given a pair of operators $T$ and $L$, define an algebra of
polynomials
\begin{equation} R_{T,L}:=\{\ p(z)\in\C[z]\ |\ Tp(L)\subset {\cal D} T\
\}\ .
\end{equation} 

Now assume Equation~\eqref{bispeqn}.  For reasons explained in
\cite{W}, it suffices to assume that $L$ and $Q$ are {\it normalized},
by which we mean that the top two coefficients are, respectively,
constant and zero.  It then follows that $L$ and $Q$ have coefficients
in the field of rational functions $\C(x)$ (resp. $z$), and that
$f(z)$ and $g(x)$ are polynomials.  For $T\in\cal A_{\psi}$, let
$a(x)$ (resp.  $\four{a}(z)$) be the leading coefficient of $T$ (resp.
$\four{T}(z,\d_z)$) and
\begin{equation}
\tilde T:=\frac 1{a(x)} T \qquad \tilde \four{T}:=\frac{1}{\four{a}(z)}\four{T}(z,\d_z).
\end{equation}  Then we immediately have the
following proposition.
\stateproposition[prop:maybe-bisp]{Given $p_1\in R_{T,L}$ and 
$p_2\in R_{\four{T},Q}$, let
\begin{equation} P_1:=\tilde{T}p_1(L)\tilde{T}^{-1}
\end{equation}
\begin{equation} P_2:=\tilde{\four{T}}p_2(Q)\tilde{\four{T}}^{-1}
\end{equation} 
and let
\begin{equation}
\phi(x,z):=\frac 1{a(x)\four{a}(z)}T(\psi(x,z))\ .
\end{equation} Then  $P_1$ and $P_2$ are normalized ordinary differential
operators, satisfying
\begin{equation}
P_1\phi=p_1(f(z))\phi\end{equation}\begin{equation}
P_{2}(z,\d_z)\phi=p_2(g(x))\phi\ .\end{equation}}

\statecorollary{Given a bispectral triple $(L,Q,\psi)$ satisfying
Equation~\ref{bispeqn}, then for $T\in {\cal A}_{\psi}$ such that the
rings $R_{T,L}$ and $R_{\four{T},Q}$ are non-trivial, the ring
\begin{equation}
{\cal R}_{T,L}:=\{ \tilde{T}p(L)
\tilde{T}^{-1}\  |\  p\in R_{T,L}\ \}
\end{equation}
 is an algebra of normalized bispectral operators, isomorphic to
$R_{T,L}$.  }

Thus one has a strategy for obtaining new bispectral algebras from a
given bispectral operator, namely to find $T\in\cal A_{\psi}$ such
that the rings $R_{T,L}$ and $R_{\four{T},Q}$ are non-trivial.  Each
such $T$ then gives a {\it bispectral Darboux transformation\/} of the
bispectral operator $p(L)$ for $p\in R_{T,L}$.

For instance, Wilson's result, though not stated as such, is
essentially the following.  Taking $\psi(x,z)=e^{xz}$, ${\cal
A}_{\psi}$ is the Weyl algebra of differential operators with
polynomial coefficients.  Indeed, in that case the map $T\mapsto
\four{T}$ is the standard antiautomorphism of the Weyl algebra
exchanging $\d$ and $x$.
\stateproposition[Wilson revisited]{Let $\lambda_1,...,\lambda_n$
be (not necessarily distinct) complex numbers and let
$p_1(x),...,p_n(x)$ be polynomials.  Let $\phi_i=
p_i(x)e^{\lambda_ix}$.  Let $\bar{K}$ be the operator
\begin{equation} \bar K(u):=e^{-x \sum\lambda_i}|Wr(\phi_1,...,\phi_n,u)|
\ ,
\end{equation}
where $|Wr(\phi_1,...,\phi_n,\cdot)|$ denotes the Wronskian operator.
(Clearly $\bar{K}$ has polynomial coefficients.)  Then
\begin{enumerate}\item The rings 
${\cal R}_{\bar{K},\d}$ and ${\cal R}_{\four{\bar{K}},\d}$ are both
nontrivial.
\item Every maximal, normalized rank one bispectral algebra is
obtained in this way.\end{enumerate}}

The purpose of the present paper is to extend part of the proposition
above to the higher rank case.  What we are able to prove at the
present time is that if we replace $\d$ by a generalized Airy operator
\begin{equation} L_0:= \d^r-a_{r-2}\d^{r-2}\cdots -a_1\d-x
\end{equation} for $r>1$ and $a_i\in\C$, then the analogue of part
one of the proposition continues to hold, at least when
$\lambda_1,...,\lambda_n$ are distinct and $p_1,...,p_n$ are
polynomials of degree one.  This analogue is the following.  The
operator $L_0$ is clearly bispectral.  Indeed, if we define
\begin{equation}
\hat f(x,z):=f(x+z)\ ,\label{eq:hatf}
\end{equation} then any $f\in \ker(L_0)$ satisfies the equations
\begin{equation}\label{eqn:eigenvector} L_0(x,\d)\hat f(x,z)=z\hat
f(x,z)\qquad L_0(z,\d_z)\hat f(x,z)=x\hat f(x,z).
\end{equation}
Note that for nonzero $f$ belonging to $\ker(L_0)$, ${\cal
A}_{\hat{f}}$ contains $\d$, since $\d(\hat{f})=\d_z(\hat{f})$, and it
contains $x$ by equation
\eqref{eqn:eigenvector}.   Thus ${\cal A}_{\hat {f}}$ is again the
Weyl algebra and so $L_0$ determines an involution of the Weyl
algebra.  Explicitly, this involution is the unique anti-automorphism
given by
\begin{equation}
\four{x}=L_0\qquad \four{\d}=\d
\ .\label{eq:four}
\end{equation}
Throughout the remainder of the paper we will consider $L_0$ fixed and
define the involution $\flat$ by Equation~\ref{eq:four}.

The main result of the paper, to be developed and proved below is

\statetheorem[-]{
Let $f_1(x),...,f_r(x)$ be a basis for the kernel of $L_0$.  Define
$\bar{K}$ to be the Wronskian operator of the $rn$ functions
\begin{equation}\label{phiij}
\phi_{i,j}:=p_i(\d_z)(\hat{f}_j)|_{z=\lambda_i}\ .
\end{equation}
 For $\lambda_1,...,\lambda_n$ distinct and $p_1,...,p_n$ polynomials
of degree one, the rings ${\cal R}_{\bar{K},L_0}$ and ${\cal
R}_{\bar{\four{K}},L_0}$ are both nontrivial and therefore
bispectral.}

One of our purposes here is to draw attention to the rather
interesting subset of the Weyl algebra consisting of those elements
$T$ for which both ${\cal R}_{T,L_0}$ and ${\cal R}_{\four{T},L_0}$
are nontrivial.  We conjecture that it corresponds by the construction
above to the set of homogeneous finite dimensional subspaces of the
space of finitely supported distributions in the complex plane, as in
Wilson's case.

By a finitely supported distribution we mean a finite linear
combination of operators of the form $\delta_{\lambda}\circ\d_z^j$,
$\lambda\in\C$, where $\delta_{\lambda} $ is the $\delta$-function
evaluating its argument at $z=\lambda$.  We call such a distribution
{\it homogeneous} if it is supported at one point, and we call a space
of distributions homogeneous if it has a homogeneous basis.

Let ${\cal S}$ be the space of finitely supported distributions in the
$z$-plane, and let $Gr_h(\cal S)$ denote the set of finite dimensional
homogeneous subspaces of ${\cal S}$.  $Gr_h(\cal S)$ appears to be a
rather interesting space from an algebro-geometric point of view.
Proposition \ref{Wilson revisited} implies that Wilson's bispectral
involution takes place on a certain quotient of $Gr_h(\cal S)$.
Notice that if equation
\eqref{bispeqn} holds, then $L$ itself belongs to ${\cal A}_{\psi}$.
Moreover, if $T$ is replaced by $Tq(L)$, where $q(z)$ is any nonzero
polynomial, then the rings $R_{T,L}$ and $R_{\four{T},L}$ remain
unchanged, as do their script counterparts.  Now it is not hard to
show that if $C\in Gr_h(\cal S)$, and $\lambda$ is a complex number
such that $\delta_{\lambda}\not\in C$ then
\begin{equation}
\bar{K}_{C\oplus\langle\delta_{\lambda}\rangle}=
\bar{K}_{C (z-\lambda)}\circ(\d-\lambda) \ ,\end{equation}
where $\bar{K}$ now refers to the operator in Proposition \ref{Wilson
revisited}.  Thus if $C\in Gr_h(\cal S)$ is a space of distributions
containing a delta function, we can use the action of $\C[z]$ on
$Gr_h(\cal S)$ to replace point $C$ with a point $C'$ that contains no
$\delta$ functions and produces the same bispectral algebras.  Let us
call $C$ {\it minimal} if it contains no $\delta$ functions.  Then
every point $C$ has a minimal representative in the above sense, so we
can think of the minimal points as either a quotient or a subset of
$Gr_h(\cal S)$.  Proposition \ref{Wilson revisited} implies that the
rank one bispectral involution may be regarded as an involution on the
minimal points of $Gr_h(\cal S)$.  Moreover, the involution takes
place on finite dimensional pieces of $Gr_h(\cal S)$, cut out by
placing an upper bound on both the order $\bar{K}$ and the degree of
its leading coefficient in Proposition \ref{Wilson revisited}.  We
conjecture that precisely the same sort of phenomenon occurs in the
generalized Airy case, though here we have only established this fact
for elements of $Gr_h(\cal S)$ having a certain form.  The natural
approach to extending these higher rank involutions to all of
$Gr_h(\cal S)$ seems to consist of first making a study of $Gr_h(\cal
S)$ as an infinite dimensional algebraic variety and then applying a
continuity argument.  This we hope to do in a future work.

\subsection{Notation}

The notation ``$:=$'' will be used to indicate an initial definition
of the object on the left hand side. Angular brackets $\langle c_i
\rangle$ indicate the linear space spanned over $\C$ by the basis
elements $c_i$. The determinant of a square matrix $M$ is denoted
$|M|$. The Wronskian matrix $Wr(\vec{v})$ of an $n$ vector $\vec v$ is
the $n\times n$ matrix whose first row is $\vec v$ and so that each
row is the derivative of the previous row. The symbols $\d_t$ will be
used to indicate $\frac{\d}{\d t}$ and $\d$ is just $\d_x$. An
ordinary differential operator $L$ in the variable $x$ is a polynomial
in the symbol $\d$ with coefficients which are functions of the
variable $x$. Often, a differential operator will be indicated by the
notation $L(x,\d_x)$ to indicate an operator in the variable $x$ and
$L(z,\d_z)$ to indicate the same operator following the change of
variables $x\to z$.

\section{The Operators $\bar K$ and $\bar \four{K}$}\label{sec:const}

First we prove a general lemma which will be used at several points
throughout the paper. Note that our distributions are assumed to act
in the $z$-variable, so that if $\psi$ is a function of $x$ and $z$,
$c(\psi)$ is a function of $x$.

\begin{lemma}[independence]{
Let $L(x,\d)$ be an $r$th order differential operator with
coefficients analytic in a neighborhood ${\cal U}\subset\C$. Let
${\cal V}\subset \C$ be an open subset and let
$\psi_1(x,z),\dots,\psi_r(x,z)$ be functions analytic in $\cal
U\times\cal V$, such that
\begin{enumerate}
\item $L(\psi_i)=z\psi_i$ for $i=1,\dots r$. \item For all
$\lambda\in\cal V$,\ $\psi_1(x,\lambda),\dots,\psi_r(x,\lambda)$ are
linearly independent functions of $x$.
\end{enumerate} Let ${\cal O}(\cal U)$ be the space of analytic
functions on ${\cal U}$ and ${\cal S}(\cal V)\subset{\cal S}$ be the
space of distributions with support in ${\cal V}$ . Then the map
\begin{eqnarray}
{\cal S}^r({\cal V})&\longrightarrow& {\cal O}(\cal U)\nonumber\\
(c_1,\dots,c_r)&\mapsto& \sum\limits_i c_i(\psi_i) \end{eqnarray} is
injective.}  Arguing by contradiction, assume $\sum\limits_i
c_i(\psi_i)=0$, with at least one $c_i$ non-zero. Let
$\lambda_1,\dots,\lambda_n$ be the points at which the $c$'s are
supported. For $j=1,\dots,n$, let $\mu_j$ be the highest order
derivative occurring among the $c$'s at $\lambda_j$. Then let $m(z)$
be the polynomial
\begin{equation}
m(z):=(z-\lambda_1)^{\mu_1}\prod_{i=2}^n(z-\lambda_i)^{\mu_i+1}.
\end{equation} Let $c_i':=c_i\circ m(z)$. Then
$c_i'=\alpha_i\delta_{\lambda_1}$ with $\alpha_i\not=0$ if and only if
$c_i$ has order $\mu_1$ at $\lambda_1$. In particular,
$c_1',\dots,c_r'$ are not all zero. However, we have
\begin{eqnarray} 0&=&m(L)(\sum_ic_i(\psi_i))\nonumber\\ &=&\sum_i c_i\circ
m(L)(\psi_i)\nonumber\\ &=&\sum_i c_i\circ m(z)(\psi_i)\nonumber\\
&=&\sum_i
\alpha_i\psi_i(x,\lambda_1)\end{eqnarray} This implies that the
$\alpha$'s are zero, which is a contradiction.
\end{lemma}

Fix a basis $\{f_j|1\leq j\leq r\}$ for $\ker L_0$. Since the
$(r-1)^{st}$ coefficient of $L_0$ is $0$, the Wronskian determinant,
$|Wr(f_1,\ldots,f_r)|$ is a non-zero constant. Assume the basis has
been chosen so that
\begin{equation} |Wr(f_1,\ldots,f_r)|=1\label{eq:wron1}.
\end{equation} 
\begin{Def}\label{def:BD} Define
$\BD\subset Gr_h(\cal S)$ to be set of all finite dimensional
subspaces $C\subset {\cal S}$ having a basis of the form
\begin{equation}
\{ c_i:=\delta_{\lambda_i}\circ(\d_z+\gamma_i)|\lambda_i\hbox{
distinct}\}. \end{equation}\end{Def}
\note The elements of $\BD$ should be thought of as identifying a
singular rational spectral curve with simple cusps at the points
$\lambda_i$ and a line bundle given by $(\gamma_1,\ldots,\gamma_n)$ in
the generalized Jacobian of that curve.

Fix an element $C\in\BD$.  For convenience, we will order the basis
for $C$, and define
\begin{eqnarray}\label{eqn:phi}
\vec\phi&:=&(c_1(\hat{f}_1),\dots,c_1(\hat{f}_r),\dots,c_n(\hat{f}_1),\dots,
c_n(\hat{f}_r))\nonumber\\ &:=& (\phi_1(x),\ldots,\phi_N(x))
\end{eqnarray} where $N:=rn$ and $\hat f(x,z)=f(x+z)$. As an immediate corollary of Lemma
\ref{independence} we have
 
\stateproposition[prop:Ndim]{
The $N$ functions, $\phi_s$, are linearly independent.}  The operator
$\bar{K}$ in the introduction is defined only up to a constant, as it
depends on a choice of basis.  To be precise we will define
\begin{equation}\label{eqn:kbar}
\bar{K}u(x)=\frac{\pm 1}{\prod_{i<j}(\lambda_i-\lambda_j)}
|Wr(\vec{\phi},u)|\ ,\label{eqn:defKbar} \end{equation} where the sign
will be specified in Theorem~\ref{asymptotics} below.

It follows from Proposition~\ref{prop:Ndim} that
$|Wr(\vec\phi)|\not=0$, so that $\bar K$ is indeed a differential
operator of order $N$. We may thus consider the ordinary differential
operator $K$ \begin{equation}
Ku(x)=\frac{|Wr(\phi_1,\ldots,\phi_N,u)|}{|Wr(\phi_1,\ldots,\phi_N)|}.
\end{equation} Note that this is the unique, monic operator of order
$N$ having the space $\langle\phi_s\rangle$ as its kernel.

The main result of this section is

\statetheorem[asymptotics]{
\begin{enumerate}
\item $\bar K$ has polynomial coefficients. \item We can (and do)
choose the sign of $\bar K$ so that the coefficient of $\d^N$ is a
monic polynomial of degree $n$. \item $\four{\bar{K}}$ also has degree
$N$.\end{enumerate}}

 Before proving this theorem, it is convenient to modify the
definition of the Wronskian matrix to take account of the present
situation. Define a sequence of operators
\begin{equation} P_j=\d^lL_0^k,\,\hbox{ where }j=rk+l,\ 0\leq l <r.
\end{equation} Then $P_j$ is monic of degree $j$. Given any vector
of functions $\vec v=(v_1(x),\ldots,v_m(x))$, define the {\it
modified\/} Wronskian matrix
\begin{equation}
\tilde{Wr}(v_1,\ldots,v_m):=\left(\matrix{\vec v\cr P_1\vec v\cr
\vdots\cr P_{m-1} \vec v}\right).
\end{equation} Since $\{P_j\}$ is related to $\{\d^j\}$ by a
unipotent transformation, we clearly have
\begin{equation} |\tilde{Wr}(v_1,\ldots,v_m)|=|Wr(v_1,\ldots,v_m)|.
\end{equation}

\begin{lemma}[some-lemma]{Given
$c=\delta_{\lambda}\circ(a\d_z+b)\in{\cal S}$, \begin{equation}
Wr(c(\hat f_1),\ldots,c(\hat f_r))=(bI+aB(x+\lambda))\Omega(x+\lambda)
\end{equation} where $\Omega(x):=Wr(f_1,\ldots,f_r)$ and \begin{equation}
B(x)=\left(\matrix{0&1&\cr &0&1&\cr &&\ddots&\ddots&\cr &&&&1&\cr
&&&&0&1\cr x&a_1&a_2&\cdots&a_{r-2}&0}\right).
\end{equation}}
 Let $v = (b, a, \underbrace{0,\dots, 0}_{r-2})$. Then
\begin{equation} (c(\hat{f}_1),
\dots, c(\hat{f}_r)) = v \Omega(x+\lambda)\ . \end{equation} Now $\Omega'(x)=
B(x)\Omega(x)$. Moreover, for $j \le r - 2$, $vB(x)^j =
(\underbrace{0,
\dots, 0}_j, b, a, 0, \dots, 0)$. Thus
\begin{equation} Wr(c(\hat{f}_1), \dots, c(\hat{f}_r)) =
\left(\matrix{ v\cr vB(x +
\lambda)\cr \vdots \cr vB(x + \lambda)^{r-1} }\right)\Omega(x+\lambda)\ .
\end{equation} But
$\left(\matrix{ v\cr vB\cr \vdots \cr vB^{r-1}}\right) = bI + aB$, so
the result follows.
\end{lemma}

\noindent Now we give the

\note[Proof of Theorem \ref{asymptotics}] For
$i=1,\dots,n$ and $j=0,\dots n-1$, let $A_{ij}$ be the $r\times r$
matrix \begin{equation} A_{ij}:=
\alpha_{ij}I+\lambda_i^r B(x+\lambda_i)\ , \label{eq:Adefa}\end{equation} where
\begin{equation}
\alpha_{ij}:=\gamma_i\lambda_i^j + j \lambda_i^{j-1}\ .
\end{equation} Given any $f\in \ker(L_0)$,
\begin{eqnarray}
(c(\hat{f}),P_1(c(\hat{f})),\dots,P_N(c(\hat{f})))&=&
(c(\hat{f}),\dots,c(\widehat{f^{(r-1)}})),\nonumber\\
&&c(z\hat{f}),\dots,c(z\widehat{f^{(r-1)}}),\dots,\nonumber\\
&&c(z^{n-1}\hat{f}),\dots,c(z^{n-1}\widehat{f^{(r-1)}}),\nonumber\\&&c(z^n\hat{f}))\
. \end{eqnarray} It then follows from Lemma \ref{some-lemma} that
\begin{equation}\label{kbar-unmasked}
\bar K u = \frac{\pm1}{\prod_{i<j}(\lambda_i-\lambda_j)} \left|
\matrix{A&\vec{U}\cr -\Lambda& P_N(u)}\right|
\end{equation} where
\begin{equation}
\vec{U}=\left(\matrix{u\cr P_1(u)\cr
\vdots\cr P_{N-1}(u)}\right),
\end{equation} {\small
\begin{equation}
\Lambda := (-\gamma_0 \lambda^n_0 -n \lambda^{n-1}_0, -\lambda^n_0,
\underbrace{0,\dots, 0}_{r-2}, \dots, -\gamma_{n-1}\lambda^n_{n-1} -
n \lambda^{n-1}_{n-1}, -\lambda^n_{n-1}, \underbrace{0, \dots,
0}_{r-2}) \end{equation}} and
\begin{equation} A:= \left(\matrix{ A_{10} & A_{20} & \dots &
A_{n0}\cr A_{11} &&&\vdots\cr
\vdots &&&\vdots\cr A_{1, n-1} &\dots & \dots & A_{n, n-1}}\right)\ 
\label{eq:Adefb}
.\end{equation} It follows immediately that $\bar{K}$ has polynomial
coefficients. Moreover, if we let $A'$ be the $n\times n$ submatrix of
$A$ in which $x$ occurs, then
\begin{equation} |A'|=\left(\prod_{i<j}(\lambda_i-\lambda_j)\right)\left(\prod_i(x+\lambda_i)\right). \end{equation}
 Letting $A''$ be the complementary $(r-1)n\times (r-1)n$ submatrix,
we have
\begin{equation} |A''|=\left(\prod_{i<j}(\lambda_i-\lambda_j)\right)^{r-1}.
\end{equation} Thus
\begin{equation} |A|=\left(\prod_{i<j}(\lambda_i-\lambda_j)\right)^rx^n+O(x^{n-1})\ ,\label{eqn:detA}
\end{equation} which proves assertion 2. To prove 3, let $K_{i,j}$
be the coefficient of $P_{ir+j}(u)$ in the determinant
Equation~\ref{kbar-unmasked}. Then
\begin{equation}
\four{\bar K}=\sum x^i\d^jK_{i,j}(L_0)\ . \end{equation} Since $\deg
K_{n,0}=\deg |A|=n$ by Equation~\ref{eqn:detA}, we must show that
\begin{equation}
\deg(K_{i,0})\le n\ \ \ \hbox{and}\end{equation}\begin{equation}
\deg(K_{i,j})< n\ \ \ \hbox{for}\ \ j>0\ . \end{equation} Labeling
the rows of $A$ as $A_0,\dots A_{N-1}$, \begin{equation} K_{ij}=\pm
|\tilde{A}_{(m)}|\ ,\end{equation} where $m=ir+j$ and
$\tilde{A}_{(m)}$ is the matrix obtained from $A$ by replacing $A_m$
with $\Lambda$. Then the same argument that proves assertion 2 proves
that $\deg(K_{i,0})\le n$. It also proves that $\deg(K_{i,j})< n$ if
$0<j<r-1$. Indeed, in that case the submatrix of $A$ which gave us the
coefficient of $x^n$ has now acquired a repeated row. But if $j=r-1$,
then we have removed one of the rows in which $x$ occurred, so in that
case the degree also drops. \endproof

\subsection{Examples} It follows from the proof of Theorem
\ref{asymptotics} that one may easily compute $\bar K$ and
$\four{\bar{K}}$. In the last column of determinant
Equation~\ref{kbar-unmasked}, replace $P_{rj+k}(u)$ with $z^j
\xi^k$. The determinant is then a polynomial in $x,z,\xi$. To get
$\bar K$, we make the replacement $z^j \xi^k\mapsto \d^k L_0^j$, and
to get $\four{\bar{K}}$ we make the same replacement after first
switching $x$ and $z$.
\note[Case $n=1$ and $r=2$]
In this case, we have $c=\delta_{\lambda}\circ(\d+\gamma)$ and
$L_0=\d^2-x$. The determinant is
\begin{equation} -\left|\matrix{
\gamma&1&1\cr \lambda + x&\gamma&
\xi\cr 1 + \gamma\lambda&\lambda&z}\right|=
-\xi+\gamma+\gamma^2\lambda-\lambda^2-\lambda x -(\gamma^2-\lambda)z
+xz \end{equation} This gives
\begin{eqnarray}
\bar{K}&=&(x+\lambda-\gamma^2)\d^2-\d+(\gamma^2-\lambda-x)x+\gamma+
\gamma^2\lambda -\lambda^2-\lambda x\ ,\label{eq:Kbarex} \\
\four{\bar{K}}&=&(x-\lambda)\d^2 -\d
+(\lambda-x)x+\gamma+\gamma^2\lambda- \lambda^2-(\gamma^2-\lambda) x\
.\label{eq:Kflatex} \end{eqnarray} The point, which we will
generalize, is that $\four{\bar{K}}$ is obtained from $\bar{K}$ by the
involutive map \begin{equation} (\lambda,\gamma)\mapsto
(\gamma^2-\lambda,\gamma)\ .\label{eq:mapex} \end{equation} In other
words, $\four{\bar{K}_C}=\bar{K}_{C^{\beta}}$ for some $C^{\beta}$.

\note[Case $n=1$ and $r=3$]
When $r=3$, the vacuum has a parameter,
\begin{equation}
L_0=\partial^3-a \partial -x\ .
\end{equation}
With $c$ as in the previous example, the relevant determinant is
\begin{eqnarray}
 \left|\matrix{ \gamma & 1 & 0 & 1 \cr 0 & \gamma & 1 & \xi \cr
   \lambda + x & a & \gamma & {\xi^2} \cr 1 + \gamma\,\lambda &
   \lambda & 0 & z \cr }\right|&=& a - {{\gamma}^2} +
   a\,\gamma\,\lambda - {{\gamma}^3}\,\lambda - {{\lambda}^2} -
   \lambda\,x +\nonumber\\ & & \gamma\,\xi - {{\xi }^2} - a\,\gamma\,z
   + {{\gamma}^3}\,z + \lambda\,z + x\,z
\end{eqnarray}
This gives
\begin{eqnarray}
\bar{K}&=&
 (x - a \gamma + \gamma^3 + \lambda ) \partial^3 - \partial^2 + (
     \gamma + a^2 \gamma - a\gamma^3 - a \lambda - a x ) \partial
     \nonumber\\ &&+a - \gamma^2 + a \gamma \lambda - \gamma^3 \lambda
     - \lambda^2+ a \gamma x - \gamma^3 x - 2 \lambda x - x^2
\nonumber\\
\four{\bar{K}}&=&
(x -\lambda )\partial^3 - \partial^2 + ( \gamma + a \lambda - a x )
 \partial + a - \gamma^2 + a\,\gamma\,\lambda -\nonumber\\&&
 \gamma^3\,\lambda - \lambda^2 - a\,\gamma\,x + \gamma^3\,x +
 2\,\lambda\,x - x^2
\end{eqnarray}  These correspond under the involution
\begin{equation} (\lambda,\gamma)\mapsto
(a\gamma - \gamma^3 - \lambda,\gamma)\ .
\end{equation}

\section{Bispectrality}

Our aim now is to prove that for $\bar K=\bar K_C$ as defined in
Equation~\ref{eqn:kbar}, the rings $R_{\bar{K},L_0}$ and
$R_{\four{\bar{K}},L_0}$ are nontrivial.  First consider the stablizer
of $C$ in the polynomial ring
\begin{equation} R_C:=\left\{p(z)\in\C[z]|c\circ p\in C\ \hbox{for all}\
c\in C\right\}. \end{equation} It follows from the formula
\begin{equation}
\delta_{\lambda}\circ\d_z\circ p(z)=p(\lambda)\circ
\delta_{\lambda}\circ\d_z+p'(\lambda)\circ\delta_{\lambda}
\end{equation}
that \begin{equation}\label{rsubc} R_C=\{p(z)|p'(\lambda_i)=0,\ 1\leq
i \leq n\}.\label{eq:alt.def}
\end{equation}
\begin{lemma}[lem:get-stable]{
$R_{\bar{K},L_0}=R_C$.}

It is a standard fact that the left ideal ${\cal D} K$ is precisely
the set of operators which annihilate the kernel of $K$.  In
particular, $p(z)$ belongs to $R_{\bar{K},L_0}$ if and only if $K
p(L_0)(c(\hat{f}))=0$ for every $f\in \ker(L_0)$ and $c\in C$.  Since
$K p(L_0)(c(\hat{f}))=\bar{K}c(p(z)\hat{f})$, this is the same has
saying $C\circ p(z)(\widehat{ker(L_0)})\subset C(\widehat{ker(L_0)})$.
By Lemma
\ref{independence}, this occurs if and only if $C\circ p(z)
\subset C$.\end{lemma}

It now follows from \eqref{rsubc} that $R_{\bar{K},L_0}\not=\C$ and in
fact contains elements of every sufficiently high degree.
Consequently,
\begin{equation}
{\cal R}_{\bar K,L_0}:=\left\{Kp(L_0)K^{-1}|p\in R_{\bar
K,L_0}\right\}
\label{eq:calR}
\end{equation}
is a commutative ring of ordinary differential operators of rank $r$.
We remark that this ring is the iterated Darboux transformation
\cite{Matveev} of $\C[L_0]$ by the zero eigenfunctions $\phi_s$.

In fact, it is well-known that for any point $C\in Gr_h(\cal S)$,
$R_C$ is non-trivial.  Thus, to prove that $R_{\four{\bar{K}_C},L_0}$
is also non-trivial it suffices to find a space $C^{\beta}\in
Gr_h(\cal S)$ such that $\four{\bar{K}_C}=\bar{K}_{C^{\beta}}$.

\begin{Def}\label{def:beta} Let $C^{\beta}\subset{\cal S}$ be the
subspace $$ C^{\beta}:=\{c\in{\cal S}|c(\hat f)\in\ker \four{\bar K}
\ \forall f\in\ker L_0\}. $$ \end{Def} An important alternative
definition of $C^{\beta}$ is given as follows.

\begin{lemma}{
$C^{\beta}=\{c\in{\cal S}|c\circ{\bar K(z,\d_z)}=0 \}$.}  Using the
definition \eqref{pair} of the involution $\flat$, $c$ belongs to
$C^{\beta}$ if and only if \begin{eqnarray} 0&=&\four{\bar K}\circ
c(\hat f)\nonumber\\ &=&c\circ\four{\bar K}(\hat f)\nonumber\\
&=&c\circ \bar K(z,\d_z)(\hat f)\ .
\end{eqnarray} By Lemma \ref{independence}, this occurs if and only
if $c\circ \bar K(z,\d_z)=0$.
\end{lemma}

\begin{lemma}{ $C^{\beta}$ is at most $n$-dimensional.}
 By Theorem~\ref{asymptotics}~(3), $\four{\bar K}$ has degree
$N=rn$. Thus, its kernel is at most $N$ dimensional. However, by
Proposition~\ref{prop:Ndim}, the vector space $\{c(\hat f)|c\in
C^{\beta}\ f\in\ker L_0\}$ has dimension $r\cdot \dim C^{\beta}$.
Since this is contained in the kernel of $\four{\bar K}$ by
definition, $\dim C^{\beta}\leq n$.
\end{lemma}

Now we make a general observation about Wronskians, which will enable
us to conclude that the dimension of $C^{\beta}$ is exactly $n$.

\begin{theorem}{ Let $\vec{g}:=(g_1(z),\dots,g_m(z))$ be a vector of
analytic functions of $z$. Let $w(z):=|Wr(g_1(z),\dots,g_m(z))|$ and
let $Q(z,\d_z)$ be the differential operator whose application to an
arbitrary function $u(z)$ is given by
$Qu=|Wr(g_1(z),\dots,g_m(z),u)|$.  Given $\lambda\in\C$ such that
$w(\lambda)=0$ and $w'(\lambda)\not=0$, there exists a unique
distribution of the form $c=\delta_{\lambda}\circ(\d + \gamma)$ such
that $c\circ Q=0$.  Moreover, if we set
\begin{equation}
Q=w(z)\d_z^m-w'(z)\d_z^{m-1}+v(z)\d_z^{m-2}+O(\d_z^{m-3})\ ,
\end{equation} then
\begin{equation}\label{eqn:gamma}
\gamma=\frac{v(\lambda)-w''(\lambda)} {w'(\lambda)}\ . \end{equation}
} We have
\begin{equation} A_1:=\left(\d\circ
Q\right)(u)\bigg|_{z=\lambda}=\left| \matrix{H& \vec{u}\cr
\vec{g}(m+1) & u^{(m+1)}(\lambda)}\right| \end{equation} and
\begin{equation} A_2:= Q(u)\bigg|_{z=\lambda}=\left|\matrix{ H &
\vec{u}\cr \vec{g}(m) & u^{(m)}(\lambda)}\right|
\end{equation} where $\vec{g}(j)$ is the $1\times m$ vector
\begin{equation}
\vec{g}(j):=\frac{d^j}{dz^j}(g_1,\ldots,g_m)\bigg|_{z=\lambda},
\end{equation}
$H$ is the $m\times m$ matrix
\begin{equation} H:=\left(\matrix{\vec{g}(0)\cr \vdots\cr
\vec{g}(m-1)}\right)
\end{equation} and
$\vec{u}$ is the $m\times 1$ vector
\begin{equation}
\vec{u}:=\left(\matrix{u(\lambda)\cr u'(\lambda)\cr \vdots\cr
u^{(m-1)}(\lambda)}\right).
\end{equation} Since $|H|=0$, the ``bottom right" corner of each
matrix above does not affect the determinant and \begin{equation}
A_1+\gamma A_2=\left|\matrix{ H & \vec{u}\cr
\vec{g}(m+1)+\gamma\vec{g}(m) & 0}\right|. \end{equation} This
determinant is zero (for arbitrary $u$) if and only if
$\vec{g}(m+1)+\gamma\vec{g}(m)\in V$ for \begin{equation}
V:=\langle\vec{g}(s)|0\leq s \leq m-1\rangle\subset \C^m.
\end{equation} Furthermore, since $w'(\lambda)\not=0$, the vectors
$\vec{g}(s)$ for $s=0,\ldots,m$ span all of $\C^m$. Consequently,
there is a unique $\gamma$ such that \begin{equation}
\vec{g}(m+1)=-\gamma {g}(m)+
\sum_{i=0}^{m-2}\alpha_i\vec{g}(i)\ .\end{equation} Such a $\gamma$
gives us our $c$. On the other hand, a direct calculation shows that
the highest order terms of $c_{\lambda}\circ Q$ vanish only if
$\gamma$ has the form $(\ref{eqn:gamma})$.
\end{theorem}

In the rank one case, the bispectral involution was shown to exchange
the roles of the $\tau$-function, which is zero when the $x$-orbit
leaves the big cell in the grassmannian, and the polynomial $q(z)$
which is zero at the singular points of the spectral curve
\cite{cmbis,W}. This pair of polynomials will be seen to satisfy the
same relationship and will also play an important role in the present
paper.

\begin{Def}\label{def:tau} Let $\tau(x)=\tau_C(x)$ be the
coefficient of the $\d^N$ in $\bar K$. In particular,
\begin{equation}
\tau(x):=\frac{\pm1}{\left(\prod_{i<j}(\lambda_i-\lambda_j)\right)^r}|A|
\end{equation} is a monic polynomial of degree $n$
where $A$ is as defined by $(\ref{eq:Adefa})$ and
$(\ref{eq:Adefb})$. Let
$q_C(z)=q(z):=\prod_{i=1}^n(z-\lambda_i)$. \end{Def}

Let $\BD(n)$ be the subset of $\BD$ consisting of subspaces of
dimension $n$.  We immediately have
\statecorollary[jackpot]{If $C\in\BD(n)$ and $\tau$ has distinct
roots, then $C^{\beta}\in\BD(n)$. Moreover, $C^{\beta}$ has support at
the roots of the polynomial $\tau$.}

Now we can deduce the main results. For convenience, we denote
$\bar{K}^{\beta}:=\bar{K}_{C^{\beta}}$,
$\tau^{\beta}(x)=\tau_{C^{\beta}}$ and denote by $\BDbeta$ the subset
\begin{equation}
\BDbeta:=\left\{C\in\BD|\tau_C\ \hbox{has distinct roots}\right\}.
\end{equation}

\begin{theorem}[thm:bdbeta]{Assume $C\in\BDbeta$. Then

\begin{enumerate}
\item $\tau^{\beta}=q$
\item $\four{\bar{K}}=\bar{K}^{\beta}$.
\end{enumerate}
Consequently, one may conclude that $\beta$ is an involution on
$\BDbeta$.  } By definition, $\four{\bar{K}}$ belongs to the ideal of
operators annihilating $c(\hat{f})$ for $f\in \ker(L_0)$ and $c\in
C^{\beta}$.  By Corollary~\ref{jackpot} and Theorem~\ref{asymptotics},
both $\four{\bar{K}}$ and $\bar{K}^{\beta}$ have degree $N$. Thus
there exists a non-zero rational function $\mu\in\C(x)$ such that
\begin{equation}
\four{\bar{K}}=\mu\bar{K}^{\beta}\ .\end{equation} Observe next that
that the the coefficients of $\four{\bar{K}}$ have degree at most
$n$. This follows from the definition of $\flat$ and the fact that
$\bar K$ has order $N$. Moreover, $c\circ\four{\bar{K}}=0$ for all
$c\in C$. It follows that the leading coefficient of $\four{\bar{K}}$
vanishes on the support of $C$. This coefficient is easily seen to be
monic, whence
\begin{equation}
\four{\bar{K}}=q(x)\partial^N+ b(x)\partial^{N-1} +
O(\partial^{N-2}). \end{equation} for some polynomial $b(x)$ of degree
at most $n$. We also have
\begin{equation}
\bar{K}^{\beta}=\tau^{\beta}(x)\partial^N
-\tau'_{\beta}(x)\partial^{N-1} + O(\partial^{N-2})\ .
\end{equation} Thus both claims are proved once we show that
$\mu=1$.

Again since $c\circ\four{\bar{K}}=0$,
\begin{equation} b(\lambda)=-q'(\lambda)\end{equation} for all
$\lambda$ in the support of $C$. Thus, there exists a constant
$\xi\in\C$, such that
\begin{equation} b(x)=-q'(x) + \xi q(x)\ .\end{equation} Then
\begin{eqnarray}
\mu\tau^{\beta}&=&q\nonumber\\
\mu{\tau^{\beta}}'&=&q'-\xi q\ .
\end{eqnarray} Taking quotients and integrating, we find that there
exists $\rho\in\C^*$ such that
\begin{equation}
\tau^{\beta}(x)=\rho q(x) e^{-\xi x}\ .
\end{equation} Since $\tau^{\beta}$ and $q$ are polynomials,
$\rho=1$ and $\xi=0$.

For any $c\in C$, $c(\hat f)\in \ker \bar K$ by definition of
$K$. Furthermore, since $\tau^{\beta}=q$, it has distinct roots
implying that $C\in \BDbeta$.  Therefore, by the results above,
$\four{(\bar K^{\beta})}=\four{(\four{\bar K})}=\bar K$ and so $c(\hat
f)\in\ker\four{(\bar K^{\beta})}$. Thus $C\subset
(C^{\beta})^{\beta}$. Since both spaces have dimension $n$, they are
equal.
\end{theorem}

\begin{theorem}[thm:bispectral!]{For any $C\in \BD$ and $\bar K$ ($=\bar K_C$), the ring ${\cal R}_{\bar
K,L_0}$ is a bispectral ring of operators.}  By
Proposition~\ref{prop:maybe-bisp} it is sufficient to show that
$R_{\bar K,L_0}$ and $R_{\bar \four{K},L_0}$ both contain non-constant
polynomials.  It follows from Lemma~\ref{lem:get-stable} that $R_{\bar
K,L_0}$ is non-trivial.  One can show that if $Q_t$ is a family of
elements in the Weyl algebra depending continuously on a parameter
$t$, $Q_t$ having a fixed order in both $\d$ and $x$ for all $t$, then
the property of $R_{Q_t,L_0}$ being nontrivial is preserved under
limits with respect to $t$.  Thus we reduce to the case that $\tau$
has distinct roots.  In such a case, Theorem~\ref{thm:bdbeta} shows
that $\bar \four{K}=\bar K_{C^{\beta}}$ for $C^{\beta}\in\BDbeta$ and
so, once again, Lemma~\ref{lem:get-stable} demonstrates that $R_{\bar
\four{K},L_0}$ is non-trivial.
\end{theorem}

\begin{Def}\label{def:newf}For any $f\in \ker L_0$, let
\begin{equation} {\bf f}_C(x,z):=\frac{1}{q(z)} K \hat f(x,z).
\end{equation}\end{Def} 

The map $f\mapsto {\bf f}_C$ takes $\ker L_0$ to an $r$-dimensional
space of common eigenfunctions for the operators $Q_p=Kp(L_0)K^{-1}\in
{\cal R}_{\bar K,L_0}$, satisfying \begin{equation} Q_p {\bf
f}_C(x,z)=p(z) {\bf f}_C(x,z).\label{eqn:q_p}
\end{equation}
As with the rank one bispectral involution \cite{W}, the action of
$\beta$ exchanges the roles of the spectral and spatial parameter in
the eigenfunctions, as demonstrated by the following result.

\begin{proposition}[relate-eigens]{For any $f\in \ker L_0$ and any
$C\in \BDbeta$ the eigenfunctions ${\bf f}_C$ and ${\bf
f}_{C^{\beta}}$ are related by the formula
\begin{equation} {\bf f}_C(x,z)={\bf f}_{C^{\beta}}(z,x).
\end{equation}}
\begin{eqnarray} {\bf f}_C(x,z) &=& \frac{1}{q(z)} K \hat f(x,z)\nonumber\\
&=& \frac{1}{\tau(x)q(z)} \bar K \hat f(x,z)\nonumber\\ &=&
\frac{1}{\tau(x)q(z)} \four{\bar K}(z,\d_z) \hat f(x,z)\nonumber\\ &=&
\frac{1}{\tau(x)q(z)} q(z) K^{\beta}(z,\d_z) \hat f(x,z)\nonumber\\ &=&
\frac{1}{q^{\beta}(x)} K^{\beta}(z,\d_z) \hat f(x,z) \nonumber\\ &=& {\bf
f}_{C^{\beta}}(z,x). \end{eqnarray}
\end{proposition}

\section{True Rank}\label{sec:true}

The ring of operators ${\cal R}_{\bar K,L_0}$ clearly has rank
$r$. However, if this ring is contained in a {\it larger\/}
commutative ring of lower rank, it is said to only have ``fake'' rank
$r$ \cite{LP,PW}. For example, if some other operator of order
relatively prime to $r$ commutes with the elements of ${\cal R}_{\bar
K,L_0}$, then it actually has ``true'' rank one. Since the bispectral
operators contained in rank one bispectral rings have already been
identified by Wilson
\cite{W}, it is useful to note that the ring ${\cal R}_{\bar K,L_0}$ is not
contained in a commutative ring of rank smaller than $r$.

\begin{Def} The {\it true rank\/} of an ordinary differential
operator is the rank of its centralizer in the ring of all ordinary
differential operators. The true rank of a commutative ring of
operators is the true rank of any of its elements.\end{Def}

\begin{lemma}[rankL0]{The centralizer of $L_0$ in the ring of
ordinary differential operators is the polynomial ring $\C[L_0]$.  In
particular, $L_0$ has true rank $r$.}  Since $L_0$ is contained in the
Weyl algebra of ordinary differential operators having polynomial
coefficients, its centralizer is as well.  Since $\flat$ is an algebra
antiautomorphism, it suffices to prove the result for $\four{L_0}$.
But $\four{L_0}=x$, and it is well-known that the centralizer of $x$
in the Weyl algebra is $\C[x]$.
\end{lemma}

\begin{lemma}[samerk]{If $X=Y_1Y_2$ and $\hat X=Y_2Y_1$ then $X$
and $\hat X$ have the same true rank.}  If $Q$ is an operator
commuting with $X$, then
\begin{eqnarray} 0&=&Y_2[Q,X]Y_1\nonumber\\ &=& Y_2QY_1Y_2Y_1-Y_2Y_1Y_2QY_1\nonumber\\
&=&[Y_2QY_1,\hat{X}].\label{eq:comm}
\end{eqnarray} Let $r$ be the true rank of $\hat X$. Then, by
$(\ref{eq:comm})$, we have $ord(Y_2QY_1)\equiv0 \hbox{ mod }r$. But,
$ord(Y_2QY_1)=ord(Y_2Y_1)+ord(Q)$ and since $ord(Y_2Y_1)\equiv 0\hbox{
mod }r$ we conclude that $ord(Q)\equiv 0\hbox{ mod }r$.  Therefore,
the true rank of $\hat X$ divides the true rank of $X$.  Then, by
symmetry, the true ranks are equal.\end{lemma}

\note Lemma~\ref{samerk} generalizes a previous result \cite{LP}
demonstrating that Darboux transformations preserve true rank in the
case $r=2$.

\begin{proposition}[truerk]{The ring ${\cal R}_{\bar K,L_0}$ of ordinary
differential operators has true rank $r$.}  The operator $q^2(L_0)\in
\C[L_0]$ has true rank $r$, according to Lemma~\ref{rankL0}.
Furthermore, \begin{equation} q^2(L_0)c(\hat f)=c(q^2(z)\hat
f)=0\qquad \forall c\in C,\ f\in \ker L_0. \end{equation} Thus, in
particular, $q^2(L_0)=QK$ for some $Q\in{\cal D}$.  Then the operator
$Q_{q^2}=Kq^2(L_0)K^{-1}\in {\cal R}_{\bar K,L_0}$ is given by $KQ$
and has the same true rank as $q^2(L_0)$. \end{proposition}

\section{Discussion}

\subsection{Examples} 
The following example demonstrates the construction of bispectral
algebras of rank $r=2$ in the case $n=1$.  For $r=2$, there is only
one choice for a generalized Airy vacuum: $L_0=\d^2-x$.  Thus, we may
choose the classical Airy functions $f_1(x)=Ai(x)$ and $f_2(x)=Bi(x)$
as a basis for $\ker L_0$.  Consider the one-dimensional space
$C=\langle \delta_{0}\circ (\d_z+{\gamma})\rangle$.  A polynomial
$p(z)$ stabilizes $c$ if and only if $p'(0)=0$.  Consequently,
$R_C=\C[z^2,z^3]$. From Equation~\eqref{eq:Kbarex}
$\tau(x)=x-{\gamma}^2$ and
\begin{equation}
K=\d^2+\frac{1}{{\gamma}^2-x}\d+\frac{x^2-{\gamma}-{\gamma}^2x}{{\gamma}^2-x}.
\end{equation} Then from Definition~\ref{def:newf}
\begin{equation}{\bf
f}_C(x,z)=(1-\frac{{\gamma}}{z({\gamma}^2-x)})\hat
f(x,z)+\frac{1}{z({\gamma}^2-x)}\hat f'(x,z) \end{equation} and from
Equation~$(\ref{eq:calR})$
\begin{equation} {\cal R}_{\bar K,L_0}=\C[L_4,L_6]
\end{equation} where {\small
\begin{eqnarray*} L_4 &:=&
\d^4 + {{2 \left( -2 - {{\gamma}^4} x + 2 {{\gamma}^2} {x^2} - {x^3}
\right) }\over {{{\left( {{\gamma}^2} - x \right) }^2}}} \d^2\\ &&+
{{-2 \left( -4 + 2 {{\gamma}^3} - {{\gamma}^6} - 2 {\gamma} x + 3
{{\gamma}^4} x - 3 {{\gamma}^2} {x^2} + {x^3} \right) }\over {{{\left(
-{{\gamma}^2} + x \right) }^3}}} \d \\ &&+ {x^2} - {8\over {{{\left(
-{{\gamma}^2} + x \right) }^4}}} - {{4 {\gamma}}\over {{{\left(
-{{\gamma}^2} + x \right) }^3}}} + {2\over {-{{\gamma}^2} + x}}
\end{eqnarray*} and
\begin{eqnarray*} L_6&:=&
\d^6+ {{3 \left( -2 - {{\gamma}^4} x + 2 {{\gamma}^2} {x^2} - {x^3}
\right) }\over {{{\left( {{\gamma}^2} - x \right) }^2}}}\d^4 \\ &&+
\left(-6 + {{24}\over {{{\left( -{{\gamma}^2} + x \right) }^3}}} +
{{6 {\gamma}}\over {{{\left( -{{\gamma}^2} + x \right) }^2}}}\right)
\d^3 \\ &&+\left(3 {x^2} - {{72}\over {{{\left( -{{\gamma}^2} + x
\right) }^4}}} - {{18 {\gamma}}\over {{{\left( -{{\gamma}^2} + x
\right) }^3}}} + {{6 {{\gamma}^2}}\over {{{\left( -{{\gamma}^2} + x
\right) }^2}}} + {9\over {-{{\gamma}^2} + x}} \right)
\d^2 \\ && +\left(6 x + {{144}\over {{{\left( -{{\gamma}^2} + x
\right) }^5}}} + {{36 {\gamma}}\over {{{\left( -{{\gamma}^2} + x
\right) }^4}}} - {{12 {{\gamma}^2}}\over {{{\left( -{{\gamma}^2} + x
\right) }^3}}} - {{3 \left( 3 + 2 {{\gamma}^3} \right) }\over
{{{\left( -{{\gamma}^2} + x \right) }^2}}} - {{6 {\gamma}}\over
{-{{\gamma}^2} + x}}\right)\d \\ &&-1 - {x^3} - {{144}\over {{{\left(
-{{\gamma}^2} + x \right) }^6}}} - {{36 {\gamma}}\over {{{\left(
-{{\gamma}^2} + x \right) }^5}}} + {{12 {{\gamma}^2}}\over {{{\left(
-{{\gamma}^2} + x \right) }^4}}} \\&&+ {{3 \left( 3 + 2 {{\gamma}^3}
\right) }\over {{{\left( -{{\gamma}^2} + x \right) }^3}}} + {{3
{\gamma}}\over {{{\left( -{{\gamma}^2} + x \right) }^2}}} - {{3
{{\gamma}^2}}\over {-{{\gamma}^2} + x}}
\end{eqnarray*} }

{}From Equation~\eqref{eq:mapex} we see that $C^{\beta}=\langle
\delta_{{\gamma}^2}\circ(\d_z+{\gamma})\rangle$ and so
$\tau^{\beta}(x)=x$,
$R_{C^{\beta}}=\C[(z-{\gamma}^2)^2,(z-{\gamma}^2)^3]$ and
\begin{equation} K^{\beta}=\d^2-\frac{1}{x}\d - x+
\frac{{\gamma}}{x}- {\gamma}^2. \end{equation} Most significantly,
computing ${\bf
f}_{C^{\beta}}(x,z):=\frac{1}{z-{\gamma}^2}K^{\beta}\hat f(x,z)$ it is
easily determined that \begin{equation} {\bf f}_{C^{\beta}}(x,z)={ \bf
f}_C(z,x)
\end{equation} and so the operators in ${\cal
R}_{C^{\beta}}:=K\C[(L_0-{\gamma}^2)^2,(L_0-{\gamma}^2)^3]K^{-1}$
demonstrate the bispectrality of ${\cal R}_{\bar K,L_0}$.

In the case that ${\gamma}=0$, $C$ and $C^{\beta}$ are the same and so
this is a fixed point of the involution.  This example is discussed by
Gr\"unbaum \cite{Gr} as a solution of the KP hierarchy. For arbitrary
${\gamma}$, the time dependent KP solution corresponding to this ring
is a non-vanishing rational solution \cite{V} and is described in
\cite{higherdual}.

\subsection{Sato Grassmannian and KP Flows}

The composition of the Wilson's involution on $Gr^1$ with the KP flows
was shown to be a non-isospectral symmetry of the KP hierarchy and a
linearizing map of the Calogero-Moser Particle System \cite{cmbis}.
The higher rank involution presented here may similarly have
interesting dynamic properties.  Therefore, the embedding of
$Gr_h({\cal S})$ into the Sato grassmannian $Gr^r$ used in
constructions of rank $r$ KP solutions \cite{PW} may be of interest.

Essentially, the construction above is a special case of the ``dual
construction'' \cite{higherdual} which associates a point in $Gr^r$ to
subspaces of ${\cal S}^r$ whose dimension is divisible by $r$.  Let
$$F(x):=(Wr(f_1(x),\ldots,f_r(x)))^{-1}$$ and ${\cal F}_j(x)$ denote
the $j^{th}$ column of this matrix. Define the composition of an
$r$-vector $\vec{v}=(v_1(z),\ldots,v_r(z))$ with an element $c\in{\cal
S}$ to be the element $c\circ \vec{v}=(c\circ v_1,\ldots,c\circ
v_r)\in {\cal S}^r$:
\begin{equation} (c\circ
\vec{v})\left(\vec{w}\right):=c(\vec{v}\cdot \vec{w}).
\end{equation} Then
$\hat C\subset {\cal S}^r$ is the $N$ dimensional subspace
\begin{equation}
\hat C:=\langle c\circ {\cal F}_j(z)|c\in C\ 1\leq j\leq r\rangle.
\end{equation}

To the vector space $\hat C$, we associate the point\footnote{In the
case that $q(z)$ has a root on the unit circle $S^1$, the point $W$
achieved in this way is not contained in $Gr^r$ as formulated in
\cite{PW} since its elements were taken to be $L^2(S^1)$.  This
subtlety need not concern us here, however.} $W\in Gr^r$
\begin{equation} W:=\frac{1}{q(z)}\overline{V_{\hat C}}
\end{equation} where $V_{\hat C}$ is the null space of $\hat C$ in
$(\C[z])^r$ and the overline indicates Hilbert closure in
$H^r=L^2(S^1)^r$.
\stateproposition{The vector Baker function $\psi_W(x,z)$ and the
eigenfunctions functions ${\bf f}_j(x,z):=\frac{1}{q(z)}K\hat
f_j(x,z)$ corresponding to $C$ are related by the formula
\begin{equation}
\psi_W(x,z)=\left({\bf f}_1(x,z),\ldots,{\bf f}_r(x,z)\right)\cdot
F^{-1}(z) .
\end{equation} Furthermore, the ring of ordinary differential
operators associated to $W$ is the ring ${\cal R}_{\bar K,L_0}$.}

Then, if we denote by $W$ the point in $Gr^r$ corresponding to $C$ and
by $W^{\beta}$ the point corresponding to $C^{\beta}$,
Proposition~\ref{relate-eigens} implies

\statecorollary{The action of the involution $\beta$ 
on vector baker functions is given by the formula
\begin{equation}\psi_W(x,z)=\psi_{W^{\beta}}(z,x)\cdot
F^{-1}(x)\cdot F(z).\end{equation}}

\subsection{Bispectral Flows}

In light of the remarks of the preceding section, it follows from the
results of \bibref[PW] that composing the distributions in $C$ with
the function $e^{tz^k}$ induces the $rk^{th}$ flow of the KP
hierarchy.  In particular, ${\cal L}=K(t)L_0^{1/r}K^{-1}(t)$ satisfies
the equation $$
\frac{\d}{\d t}{\cal L}=\left[({\cal L}^{rk})_{+},{\cal L}\right].
$$

For $C\in\BDbeta$ and $t$ sufficiently small, $C\circ
e^{tz^k}\in\BDbeta$.  Consequently, one may conclude that the KP flows
whose indices are zero modulo $r$ are a local symmetry of these
bispectral solutions.  Furthermore, as in \bibref[cmbis], the
composition of the KP flows with the bispectral involution is a
non-isospectral symmetry.  Finally, it can be shown that the
bispectral involution is a symplectic map linearizing the Hamiltonian
system determining the motion of the poles of these KP solutions.  (We
intend to discuss these results further in an upcoming paper.)

\subsection{General Remarks}

Associated to any commutative ring of ordinary differential operators
of rank $r$ is its geometric spectral data, which are a rank $r$
vector bundle over a complex projective curve
\cite{BC,70yrs}. The geometric spectral data corresponding to rank
one bispectral rings are ``unicursal'' rational curves with any line
bundle \cite{W}. That is, any rational curve having singularities in
the form of cusps with any line bundle is the geometric spectral data
of some bispectral ring of operators. The spectral curves associated
to the bispectral rings constructed in this paper are still rational
curves with singularities only in the form of cusps.  However, the
bispectral rings presented here do not correspond to arbitrary rank
$r$ vector bundles, but rather only very special bundles.  Indeed, for
given $n$ and $r$, our construction depends on $n+r-2$ parameters,
which is in general fewer than the number of moduli for a bundle of
rank $r$ on a singular rational curve of genus $n$.

Previous work has indicated a relationship between bispectrality and
polynomial $\tau$-functions \cite{cmbis,Z2}. The results above
continue to support such a relationship
(Definition~\ref{def:tau}). Furthermore, the bispectral involution
involves not only an exchange of the variables $x$ and $z$, but of the
roles of the polynomials $\tau(x)$ and $q(z)$. Since the roots of
$\tau(x)$ indicate places that the orbit of the point $W\in Gr^r$
leave the big cell, and the roots of $q(z)$ indicate the singular
points of the spectral curve, the bispectral involution can be roughly
seen as an exchange of spectral curve with the $x$-orbit.

The rings ${\cal R}$ constructed in this paper contain operators
normalized so that they are monic and have zero second coefficient.  A
change in variable or conjugation by a non-constant function would
change this normalization but preserve bispectrality \cite{W}.

It has been noted previously \cite{Gr2} that in the case of the known
bispectral operators associated to the KP hierarchy, one is able to
determine a corresponding operator in $z$ whose eigenvalue is any
antiderivative of the $\tau$ function. The same can be shown to be
true in the present case. In particular, for every $p\in
R_{C^{\beta}}$, there is an operator $Q_p(z,\d_z)$ for which ${\bf
f}_C(x,z)$ is an eigenfunction with eigenvalue $p(x)$.  However, since
$q^{\beta}(z)=\tau(z)$, we have by Equation~$\ref{eq:alt.def}$ that
$R_{C^{\beta}}$ is exactly the ring of polynomials whose derivatives
are divisible by $\tau$.

Finally, in light of Wilson's result, it is natural to conjecture that
given $L_0$, one may associate a bispectral ring to any finite
dimensional homogeneous subspace $C\subset\cal S$ by defining ${\cal
R}_{\bar K,L_0}$ as above.  In particular, we conjecture that if we
define $C^{\beta}$ as in Definition
\ref{def:beta} for any $C\in Gr_h(\cal S)$, there always exists a  
nonzero $\mu\in {\blackboard C}(x)$ such that
\begin{equation}\label{also good}\four{\bar{K}}=\mu(x) \bar{K}^{\beta}\
,\end{equation} which would imply that
\begin{equation}
R_{\four{\bar{K}},L_0}=R_{\bar{K}^{\beta},L_0}=R_{C^{\beta}}\not=
\C\ .
\end{equation}
\def\di#1#2{\delta_{#1}\circ(\d+#2)}%
For instance, suppose $r=2$,$C=\langle
\di{\lambda_1}{\gamma_1},\di{\lambda_2}{\gamma_2}\rangle$,
$\lambda_1\not=\lambda_2$, but $C$ does not belong to $\BDbeta$.  Then
one may check directly that there are two possibilities for
$C^{\beta}$.  Letting $\hat{\lambda}$ be the unique (repeated) root of
$\tau$, $C^{\beta}$ is spanned by two distributions supported at
$\hat\lambda$, of orders either $(3,0)$ or $(2,1)$.  In both cases,
$\bar{K}^{\beta}$ has order $4$, as does $\four{\bar{K}}$.  This
implies that Equation~\eqref{also good} holds.

\note After the completion of this work, a preprint 
 \cite{Bakalov} was made available which presents similar results from
a Lie theoretic point of view.

\noindent{\bf Acknowledgements}\ The authors are grateful for
friendly and helpful correspondence with Emma Previato, Tasso Kaper,
Robin John Chapman and Martin Reinders.

\begin{bibliography}

\bibitem[Bakalov]{B. Bakalov, E. Horozov and M. Yakimov,
``Highest Weight Modules of $W_{1+\infty}$, Darboux transformations
and the bispectral problem'', {\tt q-alg/9601017}}

\bibitem[BC]{J.L. Burchnall and T.W. Chaundy, ``Commutative
Ordinary Differential Operators'' {\it Proc.\ London Math.\ Soc.\/},
211 (1923) pp. 420--440}

 {J.L. Burchnall and T.W. Chaundy, ``Commutative Ordinary Differential
Operators'' {\it Proc.\ Royal Soc.\ A\/}, 118 (1928) pp. 557--583}

\bibitem[DG]{J.J.
Duistermaat and F.A. Gr\"unbaum, ``Differential Equations in the
Spectral Parameter'' {\it Communications in Mathematical Physics\/}
103 (1986) pp. 177--240}

\bibitem[firstbisp]{F.A. Gr\"unbaum, ``The Limited Angle
Problem in Tomography and some Related Mathematical Problems'' in
Bifurcation theory, Mechanics and Physics, {\it D. Reidel
Publishing\/}, (1983) pp. 317--329}

\bibitem[Gr]{F.A. Gr\"unbaum, ``The Kadomtsev-Petviashvili
Equation: An Alternative Approach to the `Rank Two' Solutions of
Krichever and Novikov'' {\it Physics Letters A\/} 139 (1989) pp.
146--150}

\bibitem[Gr2]{F.A. Gr\"unbaum, ``Time-Band Limiting
and the Bispectral Problem'' {\it Communications on Pure and Applied
Mathematics\/} 157 (1994) pp. 307--328}

\bibitem[cmbis]{A. Kasman, ``Bispectral KP Solutions and
Linearization of Calogero-Moser Particle Systems'', {\it
Communications in Mathematical Physics} 172 (1995) pp. 427-448}

\bibitem[higherdual]{A. Kasman, ``Darboux Transformations from $n$-KdV
to KP'', in preparation}

{A.  Kasman, ``Rank $r$ KP Solutions with Singular Rational Spectral
Curves'', Ph.D. Thesis, Boston University (1995)}

\bibitem[LP]{G. Latham and E. Previato, ``Higher Rank Darboux
Transformations'', MSRI Preprint 05229-91, to appear in {\it Proc.\
NATO ARW Nonsingular Limits of Dispersive Waves\/}, Plenum Publishing}

\bibitem[Matveev]{V.B. Matveev, ``Darboux
transformation and explicit solutions of the Kadomtcev-Petviaschvily
Equation, Depending on Functional Parameters'' {\it Letters in
Mathematical Physics\/} 3 (1979), pp. 213-216}

\bibitem[70yrs]{E. Previato, ``Seventy Years of Spectral Curves:
1923--1993'' to appear in {\it Proceedings of CIME 1993\/}
Springer-Verlag, Lecture Notes in Physics} \bibitem[PW]{E.  Previato
and G. Wilson, ``Vector Bundles Over Curves and Solutions of the KP
Equations'' {\it Proc.\ Sympos. Pure Math.\/}, 49 (1989),
pp. 553--569}

\bibitem[SW]{G. Segal and G. Wilson, ``Loop Groups and Equations
of KdV Type'' {\it Publications Mathematiques No. 61 de l'Institut des
Hautes Etudes Scientifiques\/} (1985) pp. 5--65}

\bibitem[V]{A.P. Veselov, ``Rational Solutions of the KP Equation
and Hamiltonian Systems'' {\it Communications of the Moscow
Mathematical Society, Russian Math Surveys\/} 35:1 (1980) pp.
239--240}

\bibitem[W]{G. Wilson, ``Bispectral Commutative
Ordinary Differential Operators'', {\it J. reine angew.\ Math.\/} 442
(1993) pp. 177--204}

\bibitem[Wright]{P. Wright, ``Darboux transformations, algebraic
subvarieties of Grassmann manifolds, commuting flows and
bispectrality'', Ph.D. thesis, Berkeley 1987}

\bibitem[Z1]{J. Zubelli, ``Differential Equations in the
Spectral Parameter for Matrix Differential Operators'' {\it Physica
 D\/} 43 (1990) pp. 269--287}

\bibitem[Z2]{J. Zubelli, ``On the
Polynomial Tau Function for the KP Hierarchy and the Bispectral
Property'' {\it Letters in Mathematical Physics\/} 24 (1992) pp.
41--48}

\bibitem[ZM]{J. Zubelli and F. Magri, ``Differential
Equations in the Spectral Parameter, Darboux Transformations and a
Hierarchy of Master Symmetries for KdV'', {\it Communications in
Mathematical Physics\/} 141 (1991) pp. 329--351}

 \end{bibliography}

\end{document}